\documentclass[11pt, a4paper]{article}

\usepackage[english]{babel}

\usepackage{amsmath}
\usepackage{amssymb}
\usepackage{amsthm}
\usepackage{mathabx}
\usepackage{tikz}
\usepackage{tikz-cd}
\usepackage{todonotes}
\usepackage{hyperref}
\usepackage{cleveref}
\usepackage{footmisc}
\usepackage{bbm}
\usepackage{afterpage}
\usepackage{mathrsfs}
\usepackage{sseq}
\usepackage{extarrows}
\usepackage{enumitem}
\usepackage{stmaryrd}
\usepackage{pdflscape}

\usepackage{enumitem}

\newcommand{\mult}{\mathrm{mult}}

\newcommand{\red}{{\mathrm{red}}}

\DeclareMathOperator{\Spec}{Spec}
\DeclareMathOperator{\Spf}{Spf}
\DeclareMathOperator{\Syn}{Syn}

\newcommand{\sm}{{\mathrm{sm}}}

\newcommand{\Ell}{\mathrm{Ell}}
\newcommand{\et}{\mathrm{\acute{e}t}}

\newcommand{\MF}{\mathrm{MF}}
\newcommand{\mf}{{\mathrm{mf}}}

\newcommand{\ord}{{\mathrm{ord}}}

\newcommand{\Tate}{{\mathrm{Tate}}}

\newcommand{\Sph}{\mathbf{S}}

\DeclareMathOperator{\ku}{ku}
\DeclareMathOperator{\ko}{ko}
\DeclareMathOperator{\KU}{KU}

\DeclareMathOperator{\BP}{BP}
\DeclareMathOperator{\KO}{KO}
\DeclareMathOperator{\MU}{MU}
\DeclareMathOperator{\MString}{MString}
\DeclareMathOperator{\TMF}{TMF}

\DeclareMathOperator{\Tmf}{Tmf}
\DeclareMathOperator{\tmf}{tmf}

\newcommand{\s}{{\mathrm{s}}}

\DeclareMathOperator{\Ab}{Ab}

\newcommand{\Alg}{{\mathrm{Alg}}}

\newcommand{\BT}{{\mathrm{BT}^p}}
\newcommand{\BTn}{{\mathrm{BT}^p_n}}
\newcommand{\BTone}{{\mathrm{BT}^p_1}}
\newcommand{\BTtwo}{{\mathrm{BT}^p_2}}

\newcommand{\CAlg}{{\mathrm{CAlg}}}
\DeclareMathOperator{\CRing}{CRing}

\DeclareMathOperator{\Fun}{Fun}

\newcommand{\FG}{{\mathrm{FG}}}

\newcommand{\Mod}{\mathrm{Mod}}

\DeclareMathOperator{\fDM}{fDM}
\DeclareMathOperator{\DM}{DM}

\newcommand{\Sp}{\mathrm{Sp}}
\newcommand{\Spc}{\mathcal{S}}
\renewcommand{\top}{\mathrm{top}}

\renewcommand{\lim}{{\mathrm{lim}\,}}

\newcommand{\Hom}{\mathrm{Hom}}

\DeclareMathOperator{\Map}{Map}

\newcommand{\id}{\mathrm{id}}

\newcommand{\op}{\mathrm{op}}

\newcommand{\Z}{\mathbf{Z}}

\DeclareMathOperator{\Ext}{Ext}

\newcommand{\Tors}{{\mathfrak{T}\mathrm{ors}}}
\newcommand{\Free}{{\mathfrak{F}\mathrm{ree}}}

\renewcommand{\b}{\mathfrak{B}}
\newcommand{\C}{\mathcal{C}}

\newcommand{\E}{\mathbf{E}}
\newcommand{\EE}{\mathcal{E}}

\newcommand{\F}{\mathbf{F}}

\newcommand{\G}{\mathbf{G}}

\newcommand{\h}{\mathrm{h}}

\renewcommand{\j}{\mathrm{j}^1}

\newcommand{\K}{\mathrm{K}}

\newcommand{\M}{\mathcal{M}}

\renewcommand{\O}{\mathscr{O}}

\renewcommand{\P}{\mathcal{P}}

\newcommand{\Q}{\mathbf{Q}}
\newcommand{\QQ}{\mathcal{Q}}

\newcommand{\R}{\mathbf{R}}

\renewcommand{\s}{{\mathrm{j}^2}}

\newcommand{\T}{\mathrm{T}}
\renewcommand{\t}{{\mathrm{t}}}

\renewcommand{\u}{\mathrm{u}}

\newcommand{\U}{\mathrm{U}}

\newcommand{\x}{\mathfrak{X}}
\newcommand{\xx}{\mathsf{X}}

\newcommand{\al}{\alpha}
\newcommand{\be}{\beta}

\usepackage{mathtools}
\DeclarePairedDelimiter{\ceil}{\lceil}{\rceil}

\theoremstyle{theorem}\numberwithin{equation}{section}
\newtheorem{theorem}[equation]{Theorem}
\crefname{theorem}{{th}.\!\!}{{ths}.\!\!}
\Crefname{theorem}{{Th}.\!\!}{{Ths}.\!\!}
\newtheorem{theoremalph}{Theorem}

\crefname{theoremalph}{{th}.\!\!}{{ths}.\!\!}
\Crefname{theoremalph}{{Th}.\!\!}{{Ths}.\!\!}

\Crefname{problem}{{Prb}.\!\!}{{Prbs}.\!\!}
\newtheorem{prop}[equation]{Proposition}
\Crefname{prop}{{Pr}.\!\!}{{Prs}.\!\!}
\newtheorem{lemma}[equation]{Lemma}
\Crefname{lemma}{{Lm}.\!\!}{{Lms}.\!\!}
\newtheorem{cor}[equation]{Corollary}
\Crefname{cor}{{Cor}.\!\!}{{Cors}.\!\!}
\newtheorem{conjecture}[equation]{Conjecture}
\Crefname{conjecture}{{Conj}.\!\!}{{Conjs}.\!\!}

\theoremstyle{definition}\numberwithin{equation}{section}
\newtheorem{mydef}[equation]{Definition}
\Crefname{mydef}{{Df}.\!\!}{{Dfs}.\!\!}

\Crefname{recall}{{Rcl}.\!\!}{{Rcls}.\!\!}

\Crefname{construction}{{Con}.\!\!}{{Cons}.\!\!}

\Crefname{ass}{{As}.\!\!}{{As}.\!\!}
\newtheorem{notation}[equation]{Notation}
\Crefname{notation}{{Nt}.\!\!}{{Nts}.\!\!}

\Crefname{situation}{{St}.\!\!}{{Sts}.\!\!}

\theoremstyle{remark}\numberwithin{equation}{section}
\newtheorem{example}[equation]{Example}
\Crefname{example}{{Ex}.\!\!}{{Exs}.\!\!}

\Crefname{nonexample}{{NonEx}.\!\!}{{NonEx}.\!\!}

\Crefname{claim}{{Clm}.\!\!}{{Clms}.\!\!}
\newtheorem{remark}[equation]{Remark}
\Crefname{remark}{{Rmk}.\!\!}{{Rmks}.\!\!}
\newtheorem{question}[equation]{Question}
\Crefname{question}{{Q}.\!\!}{{Qs}.\!\!}

\Crefname{idea}{{Id}.\!\!}{{Ids}.\!\!}

\Crefname{warn}{{Warn}.\!\!}{{Warns}.\!\!}

\Crefname{figure}{{Fig.}\!\!}{{Figs.}\!\!}
\Crefname{footnote}{{Fn.}\!\!}{{Fn.}\!\!}

\Crefname{part}{{\textsection}\!\!}{{\textsection}\!\!}
\Crefname{chapter}{{\textsection}\!\!}{{\textsection}\!\!}
\Crefname{section}{{\textsection}\!\!}{{\textsection}\!\!}
\Crefname{subsection}{{\textsection}\!\!}{{\textsection}\!\!}
\Crefname{appendix}{{\textsection}\!\!}{{\textsection}\!\!}
\begin{document}
\title{Constructing and calculating Adams operations on\\ dualisable topological modular forms}
\author{Jack Morgan Davies}

\maketitle

\begin{abstract}
We construct Adams operations $\psi^k$ on the cohomology theory $\Tmf$ of dualisable topological modular forms after inverting $k$; the first such multiplicative stable operations on this cohomology theory. These Adams operations are then calculated on the homotopy groups of $\Tmf$ using a combination of descent spectral sequences and Anderson duality. Applications of these operations are then given, including constructions of connective height $2$ analogues of \emph{Adams summands} and \emph{image-of-$J$ spectra}.
\end{abstract}



\setcounter{tocdepth}{3}
\tableofcontents

\addcontentsline{toc}{section}{Introduction}
\section*{Introduction}

In this article, the extraordinary cohomology theory $\Tmf$ of \emph{topological modular forms} is equipped with \emph{Adams operations}, compatible with the classical Adams operations on topological $K$-theory (\Cref{tha,thaversiontwo}). These are the first nonidentity multiplicative operations on $\Tmf$. This is made possible using a combination of a powerful theorem of Lurie in spectral algebraic geometry and a careful application of Goerss--Hopkins obstruction theory. These Adams operations are then calculated on the homotopy groups of $\Tmf$ (\Cref{thb}). Finally, we make use of the relationship between $\Tmf$ and topological $K$-theory to construct height $2$ analogues of the \emph{connective Adams summand} and \emph{connective image-of-$J$ spectrum} and prove some basic facts about these spectra (\Cref{thd,the}). 

\subsection*{Motivation}

Adams operations are some of the most utilised power operations in homotopy theory. This is exemplified by the work of Adams counting the number of vector fields on spheres \cite{vectorfieldsonspheres}, by Adams--Atiyah giving a ``postcard-sized'' proof of the Hopf invariant one theorem \cite{adamsatiyah}, and by Quillen calculating the algebraic $K$-theory of finite fields \cite{quillenkff}. In these three examples, Adams operations arise as operators on \emph{topological $K$-theory} $\KU$. In this article, we study Adams operations on another extraordinary cohomology theory $\Tmf$, called \emph{topological modular forms}, which has received much attention in recent years. Many view $\Tmf$ as a natural higher height analogue of topological $K$-theory, due to its relationships with number theory (through its connection to modular forms) as well as differential geometry and physics (through the string orientation $\sigma\colon \MString\to \Tmf$), and its ability to help with computations in stable homotopy theory; see \cite{handbooktmf} for more details and references. \\

The construction of $\Tmf$ arises from the study of elliptic cohomology theories and generalised elliptic curves. To motivate our study of Adams operations on $\Tmf$, let us first reimagine the classical Adams operations on $\KU$ through an algebro-geometric lens.\\

One can recover complex topological $K$-theory $\KU$ using only the multiplicative group scheme ${\G}_m=\Spec\Z[ t^{\pm}]$. In fact, one can recover the multiplicative stable homotopy type, also known as the $\E_\infty$-ring, which represents the cohomology theory $\KU$. There is a moduli stack $\M_{\G_m}$ of \emph{forms of $\G_m$} upon which there exists an étale sheaf $\O^\mult$ of $\E_\infty$-rings, constructed using spectral algebraic geometry. When evaluated on the \'{e}tale open $\Spec \Z\to \M_{\G_m}$ defined by the multiplicative group $\G_m$ over $\Z$, we obtain $K$-theory $\O^\mult(\G_m/\Spec \Z)=\KU$. Moreover, the functoriality of $\O^\mult$ means that automorphisms of $\G_m$ over $\Spec \Z$ induce automorphisms of $\KU$. For example, the inversion isomorphism $[-1]\colon \G_m\to \G_m$ (defined by sending $t$ to $t^{-1}$) produces the Adams operation $\psi^{-1}\colon \KU\to \KU$ on $K$-theory; this is the familiar $C_2$-action sending a complex vector bundle to its conjugate. To obtain more Adams operations on $\KU$, say $\psi^k$ for each integer $k$, one might try to extend $\O^\mult$ to be functorial with respect to more endomorphisms of $\G_m$. This is not an unreasonable request, as the cohomology theory $\KU$ only really depends on the \emph{formal group} $\widehat{\G}_m$ associated with $\G_m$, by the classical \emph{Landweber exact functor theorem}, and many endomorphisms of $\G_m$ produce automorphisms on $\widehat{\G}_m$. Once this extended functoriality is achieved, the Adams operations\footnote{These operations are often called \emph{stable} Adams operations, as they are defined as maps of spectra. As these are the only kind of Adams operations we will consider in this article, let us forgo the adjective \emph{stable}.} $\psi^k$ on $\KU[\tfrac{1}{k}]$ can be obtained by applying $\O^\mult$ to the $k$-fold multiplication map on $\G_m$. This blueprint is carried out and discussed in detail in \cite[\textsection6.4]{luriestheorem} for $\KU$ completed at a prime.\\

Similar constructions can also be considered for $\Tmf$. Indeed, the definition of $\Tmf$ is as the global sections of the celebrated étale sheaf $\O^\top$ of $\E_\infty$-rings on the moduli stack of generalised elliptic curves $\M_\Ell$. This sheaf was originally constructed by Goerss--Hopkins--Miller \cite{bourbakigoerss} and takes values in \emph{elliptic cohomology theories}. This means that for each affine \'{e}tale map $E\colon \Spec R\to \M_\Ell$, the cohomology theory $\O^\top(E/\Spec R)$ remembers the formal group of the generalised elliptic curve $E$. The $\E_\infty$-ring $\Tmf$ is the global sections of $\O^\top$ on $\M_\Ell$, so we say $\Tmf$ is the \emph{universal elliptic cohomology theory} as it maps to all other elliptic cohomology theories given as sections of $\O^\top$. One can now ask if we can use the ``multiplication map'' on generalised elliptic curves to construct Adams operations on $\Tmf$ using $\O^\top$. Just as for topological $K$-theory, there is a construction of Adams operations on periodic topological modular forms $\TMF$ by following the above blueprint applied to the moduli stack of \emph{smooth} elliptic curves; see \cite[\textsection6.4]{luriestheorem} and \cite[\textsection2]{heckeontmf}. One might then suspect that $\Tmf$ also has Adams operations $\psi^k$ after inverting $k$, as this cohomology theory can be constructed using $\TMF$ and $\KU$. There are many subtleties to consider though, such as the lack of an honest multiplication map or group structure on generalised elliptic curves. The first goal of this article is to confirm this suspicion and show that $\Tmf$ does admit Adams operations.

\subsection*{Main results}

It is well-known that to define multiplicative maps $\psi^k$ on topological $K$-theory, one must invert $k$ (\cite[\textsection II.13]{bluebook}), as $\psi^k(u)=ku$ for the generator $u\in \pi_2\KU$, so $\psi^k(u^{-1})=\frac{1}{k}u^{-1}$. The same is true for periodic topological modular forms $\TMF$, as discussed in \cite[Th.F]{heckeontmf}. For this reason, we are content with constructing Adams operation $\psi^k$ only after inverting $k$.

\begin{theoremalph}\label{tha}
For every integer $k$, there is a morphism of $\E_\infty$-rings $\psi^k\colon \Tmf[\frac{1}{k}]\to \Tmf[\frac{1}{k}]$ and a commutative diagram of $\E_\infty$-rings
\[\begin{tikzcd}
{\Tmf[\frac{1}{k}]}\ar[r, "{\psi^k}"]\ar[d]	&	{\Tmf[\frac{1}{k}]}\ar[d]	\\
{\KU[\frac{1}{k}]}\ar[r, "{\psi^k}"]		&	{\KU[\frac{1}{k}].}
\end{tikzcd}\]
One can also replace $\KU[\frac{1}{k}]$ above with $\KO[\frac{1}{k}]$ or $\KO\llbracket q\rrbracket[\frac{1}{k}]$.
\end{theoremalph}

The morphism $\Tmf[\frac{1}{k}]\to \KU[\frac{1}{k}]$ is the \emph{evaluation at the cusp} map which on rational homotopy groups sends a modular form to the linear term in its $q$-expansion. The construction of the operations $\psi^k$ above will come from the more general $p$-complete statement.

\begin{theoremalph}\label{thaversiontwo}
For every prime $p$ and every $p$-adic unit $k\in \Z_p^\times$, in particular for every integer $k$ not divisible by $p$, there is a morphism of $\E_\infty$-rings $\psi^k\colon \Tmf_p\to \Tmf_p$ and a commutative diagram of $\E_\infty$-rings
\[\begin{tikzcd}
{\KU_p}\ar[d, "{\psi^k}"]	&	{\Tmf_p}\ar[d, "{\psi^k}"]\ar[r]\ar[l]	&	{\TMF_p}\ar[d, "{\psi^k}"]	\\
{\KU_p}			&	{\Tmf_p}\ar[l]\ar[r]				&	{\TMF_p}
\end{tikzcd}\]
where $\psi^k\colon \TMF_p\to \TMF_p$ are the Adams operations of \cite[Th.6.9 \& Df.6.16]{luriestheorem}. In particular, if $k$ is an integer not divisible by $p$, then $\psi^k\colon \Tmf_p\to \Tmf_p$ is the $p$-completion of the operation $\psi^k\colon \Tmf[\frac{1}{k}]\to \Tmf[\frac{1}{k}]$ of \Cref{tha}. Furthermore, if $p$ is odd, then restricting to the maximal finite subgroup of $\Z_p^\times$ yields an action of $\F_p^\times$ on $\Tmf_p$ such that the maps of $\E_\infty$-rings $\Tmf_p\to \KU_p$, $\Tmf_p\to \TMF_p$, and $\psi^k\colon \Tmf_p\to\Tmf_p$ are $\F_p^\times$-equivariant. One can also replace $\KU_p$ above with $\KO_p$ or $\KO\llbracket q\rrbracket_p$. 
\end{theoremalph}

These Adams operations above on $\Tmf[\frac{1}{k}]$ and $\Tmf_p$ are (to the best of the author's knowledge) the first nonidentity stable multiplicative operations on these cohomology theories. \\

These theorems come with a warning: there is no obvious compatibility between various Adams operations $\psi^k$ from \Cref{tha} or \Cref{thaversiontwo}. This means that we do not claim to have homotopies $\psi^k\circ \psi^\ell \simeq \psi^{k\ell}$, for example. Homotopies of this kind and more are further explored in \cite{realspectra} away from the prime $2$ and in \cite[Ths.C-D]{heckeontmf} away from $\Delta^{24}$, so on $\TMF$.\\

The proof of \Cref{tha,thaversiontwo} involves much more work than the construction of Adams operations on topological $K$-theory $\KU$ and periodic topological modular forms $\TMF$, as generalised elliptic curves do not always admit a multiplication map. It is only with a combination of the ideas used in the $\KU$- and $\TMF$-cases together with Goerss--Hopkins obstruction theory that we can prove \Cref{tha,thaversiontwo}.\\

As fundamental operations on an important cohomology theory, we expect the Adams operations of \Cref{tha,thaversiontwo} to become useful tools in algebraic topology. With an eye to these future applications, we proceed to calculate the effect of these operations on homotopy groups---it suffices to state the $p$-complete calculations here. 

\begin{theoremalph}\label{thb}
For every odd prime $p$, every $p$-adic unit $k\in \Z_p^\times$, and every $x\in \pi_d\Tmf_p$ with $d$ positive, we have the equality
\[\psi^k(x)=\begin{cases}	x & x\in\Tors_d \\ k^{\ceil{\frac{d}{2}}}x & x\in\Free_d \end{cases}\]
where $\Tors_d\subseteq\pi_d\Tmf_p$ is the subgroup of torsion elements and $\Free_d$ is the orthogonal subgroup of \Cref{elementsintmf}. At the prime $p=2$, the above equalities hold for all $d$ except for those positive $d$ congruent to $60$ or $156$ modulo $192$.\footnote{As is made clear in \Cref{elementsintmf} when defining $\Free_d$, there is some ambiguity at the prime $2$ and for these $d$, where an explicit basis for $\Free_d$ has not yet been found.}
\end{theoremalph}

Despite the similarity to the calculations of Adams operations on $\pi_\ast\KU$ and $\pi_\ast\KO$, the above theorem requires a much more detailed analysis.\\

With the construction and calculation of Adams operations on $\Tmf$ in hand, we start to imitate some of the classical constructions on topological $K$-theory and Adams operations now using topological modular forms. For example, one can split $p$-complete connective complex topological $K$-theory $\ku_p$ into \emph{Adams summands} $\ell=\ku^{h\F_p^\times}$ using the $\F_p^\times$-action from the $p$-adic Adams operations. One can also study the \emph{image-of-$J$} by defining a connective spectrum $\j$ as the further fixed points of Adams operations acting on $\ell$. One major advantage of the Adams operations on $\tmf$ of \Cref{tha,thaversiontwo} compared with those on periodic $\TMF$ of \cite[\textsection6.4]{luriestheorem} (and \cite[\textsection2]{heckeontmf}) is the direct comparison to the Adams operations on topological $K$-theory and the ability for one to use an $\F_p$-based Adams spectral sequence. These computations of the $\F_p$-ASS for $\j$ appear in \cite{brunerrognes} and a modified $\F_p$-ASS for $\j$ appears in \cite{syntheticj}. The analogy to $K$-theory is then used to motivate the study of connective Adams summands $\u=\tmf_p^{h\F_p^\times}$ and image-of-$J$-spectra $\s$ (again as a fibre of an Adams operation acting on $\u$) at the height $2$. To highlight the simplicity and utility of \Cref{tha,thaversiontwo,thb}, we prove the following two statements involving $\u$ and $\s$.\\

The $\E_\infty$-ring $\u$ is always a summand of $\tmf_p$, however, unlike the height 1 case, the $\u$-module $\tmf_p$ is not necessarily a sum of shifts of $\u$. Our next theorem summarises at which primes $\tmf_p$ splits into copies of $\u$.

\begin{theoremalph}\label{thd}
The inclusion of fixed points $\u\to \tmf_p$ witnesses $\tmf_p$ as a \emph{quasi-free} $\u$-module if and only if $p-1$ divides $12$. On the other hand, the map of $\E_\infty$-rings $\U\to \TMF_p$ always witnesses $\TMF_p$ as a quasi-free $\U$-module.
\end{theoremalph}

If $p-1$ does not divide $12$, so for primes $p=11$ and $p\geq 17$, we believe the next best thing to a splitting is true. More specifically, we conjecture (\Cref{secondpslittingtheorem}) that there is a cofibre sequence of $\u$-modules of the form
\[\bigoplus_{\frac{p-1}{2}} \u[?]\to \tmf_p\to \bigoplus \ell[?]\]
and provide such a cofibre sequence for the primes $p=11, 17, 19, 23$, and $37$.\\
 
Our final theorem takes advantage of the fact that the $q$-expansion map $\Tmf_p\to \KO$ commutes with our Adams operations and is $\F_p^\times$-equivariant, which allows us to explicitly compare $\s$ and $\j$.

\begin{theoremalph}\label{the}
Let $p$ be a prime. Then the unit map $\Sph_p\to \s$ detects all elements in $\pi_\ast \Sph_p$ in the $p$-primary image-of-$J$ and all elements detected by $\Sph\to \tmf_p$.
\end{theoremalph}

At the prime $p=3$, Carrick and the author have shown that the unit map $\Sph_3 \to \s$ detects more $v_2$-periodic families from $\pi_\ast \Sph_3$ than just those detected by $\tmf_3$; see \cite{heighttwojatthree}.\\

These last two theorems highlight some of the immediate applications of Adams operations on $\Tmf$, including their formal properties and computational power.

\subsection*{Outline}

The sections of this article \Cref{otopisbtextendedsection,adamsopessections,applicationssection} can be read independently, assuming the main results of the previous sections. In \Cref{otopisbtextendedsection}, we use many tools surrounding $\Tmf$ such as (spectral) algebraic geometry, elliptic cohomology theories, and Goerss--Hopkins obstruction theory; in \Cref{adamsopessections}, we use some formal computational aspects in stable homotopy theory, including some \emph{synthetic spectra}, and Anderson duality; in \Cref{applicationssection}, we use some stable homotopy theory surrounding the image-of-$J$ and elementary notions in the theory of modular forms. In some more detail:
\begin{itemize}
\item In \Cref{otopisbtextendedsection}, we prove \Cref{tha,thaversiontwo} and construct the titular Adams operations on $\Tmf$. This opens with an outline of the algebraic geometry \Cref{introddsection} used in this article. Next is a construction of Adams operations on sections of the sheaf $\O^\top$ over open substack $\M_\Ell^\sm$ of smooth elliptic curves and its complement $\M_{\G_m}$ using Lurie's theorem and spectral algebraic geometry \Cref{pdivsconstuctionssection}. In \Cref{gluingtogether}, we prove \Cref{thaversiontwo} by gluing together our operations on $\KO\llbracket q\rrbracket_p$ and $\TMF_p$ using Goerss--Hopkins obstruction theory. Finally, we prove \Cref{tha} in \Cref{newintegralstuff} by gluing together the $p$-complete Adams operations at different primes with some rational datum.
\item In \Cref{adamsopessections}, we prove \Cref{thb} and calculate our Adams operations on the homotopy groups of $\Tmf$. First, we define an explicit basis for our summands $\Free$ and $\Tors$ in \Cref{achoiceofgenerators} using the computations from \cite{brunerrognes}. Next, in \Cref{andersondualitysection} we discuss the Anderson self-duality of $\Tmf$ (as proven by Stojanoska) and the formal ramifications this self-duality implies. Finally, in \Cref{generaloperations}, we prove \Cref{thb} using some formal stable homotopy theory and the Anderson self-duality of $\Tmf$.
\item In \Cref{applicationssection}, we prove \Cref{thd,the} using the connections between $\Tmf$ and topological $K$-theory. We start \Cref{adamssummandapplication} with a proof of \Cref{thd}, which follows from our calculations of Adams operations on $\Tmf$ and some basic facts about spaces of modular forms. In \Cref{previouslydiscussedsplittingconjecturesection}, we discuss evidence for a conjecture explaining the negative cases of \Cref{thd}. In \Cref{imageofjapplication}, we prove \Cref{the}, which is a purely formal consequence of \Cref{thaversiontwo} and the classical study of the image-of-$J$.
\end{itemize}

\subsection*{Past and future work}
Operations on elliptic cohomology theories have been constructed by Baker \cite{bakerhecketwo, bakerss} and Ando \cite{andopo}, and these include Adams operations. The Adams operations in this article can be seen as global stable $\E_\infty$-versions of those previously studied. As mentioned by Baker \cite[p.6]{bakerhecketwo}, the Adams operations $\psi^k$ are determined as a multiplicative natural transformation of homology theories on $\TMF[\frac{1}{6k}]$ (which is the modern notation for classical elliptic cohomology) by the formula $\psi^k(x)=k^dx$ for $x\in \pi_{2d}\TMF[\frac{1}{6k}]$. By \Cref{thb}, we see that our operations $\psi^k$ are homotopic to those classical stable Adams operations on $\TMF[\frac{1}{6k}]$. We have also explored other operations on $\Tmf$ and related spectra. In \cite{heckeontmf}, we discuss Adams operations, Hecke operators, and Atkin--Lehner involutions on $\TMF$ as well as periodic topological modular forms with level structure. The structural results for operations on $\TMF$ are much stronger than those shown here for $\Tmf$, as we have a spectral algebro-geometric description of the former. The Adams operations on $\Tmf$ in this article are constructed from those on $\TMF$ from \cite{heckeontmf} or equivalently \cite{luriestheorem}. In \cite{realspectra}, we prove that away from the prime $2$ the Adams operations on $\Tmf$ in this article compose as expected $\psi^k\circ \psi^\ell\simeq \psi^{k\ell}$ up to homotopy. This begins to show that operations on $\Tmf$ ought to behave as those on $\TMF$, but so far these methods are rather ad hoc. As alluded to in \cite{realspectra} and implied by this article, there are also morphisms of $\E_\infty$-rings $\Tmf \to \Tmf_0(n)$ (and \textbf{not} of $\Tmf$-modules) critical to defining Hecke operations as well as connective forms of Behrens' $Q(N)$ spectra of \cite{ktwospheremark}. We will return to such constructions in future work. Finally, a further study of Adams operations on $\tmf_3$ leads to detection statements for products within the divided $\beta$-family in $\pi_\ast \Sph_3$; see \cite{heighttwojatthree}.

\subsection*{Conventions}
The language of $\infty$-categories will be used throughout, so all categorical constructions and considerations will be of the $\infty$-categorical flavour. In particular, for a scheme $X$ and a finite group acting on $X$, we will write $X/G$ for what is sometimes called the \emph{stacky quotient}. In general, we consider our algebraic geometry as occurring in the $\infty$-category $\Fun(\CRing,\Spc)$ where $\CRing$ is the $1$-category of commutative (discrete) rings. Given a prime $p$, we will also write $\widehat{\M}$ for $\M\times\Spf\Z_p$, where $\M$ is any presheaf in $\Fun(\CRing,\Spc)$. We will denote the $p$-completion of $\E_\infty$-rings with a subscript $(-)_p$. All of our discrete rings will be commutative. For an $\E_\infty$-ring $R$ and $R$-modules $M,N$, we write $F_R(M,N)$ for the internal function $R$-module in $\Mod_R$. For an integer $n$, we will write $X[n]$ for the $n$th suspension of a spectrum $X$. 

\subsection*{Acknowledgements}
I owe a debt of gratitude to my PhD supervisor Lennart Meier for all the ideas, support, and clarity he has shared with me. Thank you to Christian Carrick, Dan Isaksen, Guchuan Li, and John Rognes, who corrected (multiple times) and furthered my understanding of the homotopy groups of topological modular forms. I truly appreciate the careful reading of an anonymous referee, whose remarks and suggestions led to a better article. I have also had fruitful conversations with Gijs Heuts, Magdalena Kedziorek, Tommy Lundemo, and Yuqing Shi. Thank you also to Lennart, Theresa, and Tommy for proofreading.


\section{Constructions}\label{otopisbtextendedsection}

In this section, we will prove \Cref{tha,thaversiontwo} and construct Adams operations $\psi^k$ on $\Tmf[\frac{1}{k}]$ and $\Tmf_p$ as morphisms of $\E_\infty$-rings. In the $p$-complete case, we will use the Cartesian diagram of $\E_\infty$-rings
\[\begin{tikzcd}
{\Tmf_p}\ar[r]\ar[d]	&	{\KO\llbracket q\rrbracket_p}\ar[d]	\\
{\TMF_p}\ar[r]		&	{\KO\llparenthesis q\rrparenthesis_p}
\end{tikzcd}\]
by taking global section of \cite[Df.5.10]{hilllawson}. In particular, we will use spectral algebraic geometry to construct Adams operations on $\TMF_p$ and $\KO\llbracket q\rrbracket_p$ and then Goerss--Hopkins obstruction theory (\`{a} la Behrens \cite[\textsection12]{tmfbook} and Hill--Lawson \cite{hilllawson}) to glue these operations together on $\KO\llparenthesis q\rrparenthesis_p$. Some rational stable homotopy theory is needed at the end to patch together the various $p$-complete pieces. In fact, this sketch is an outline for this whole section:\\

In \Cref{introddsection}, we discuss the necessary algebraic geometry to define our tools. In \Cref{pdivsconstuctionssection}, we construct a version of $\O^\top$ for $\TMF_p$ and $\KO\llbracket q\rrbracket_p$ using Lurie's theorem (which originally appeared in \cite[Th.8.1.4]{taf} and is proven in \cite{luriestheorem} with extensive discussion). In \Cref{gluingtogether}, we define Adams operations on $\Tmf_p$ in an ad hoc manner by gluing together these operations on sections of $\TMF_p$ and $\KO\llbracket q\rrbracket_p$ using Goerss--Hopkins obstruction theory, thus proving \Cref{thaversiontwo}. Finally, in \Cref{newintegralstuff}, we prove \Cref{tha} using \Cref{thaversiontwo} and some rational arguments.


\subsection{Algebro-geometric background}\label{introddsection}
As our main algebro-geometric object of interest is the moduli of generalised elliptic curves, we will freely speak of stacks and Deligne--Mumford stacks; see \cite[\href{https://stacks.math.columbia.edu/tag/0ELS}{0ELS} \& \href{https://stacks.math.columbia.edu/tag/03YO}{03YO}]{stacks}. We will also need to consider formal Deligne--Mumford stacks, all of which we assume to be locally Noetherian; see \cite[\textsection8.1]{sagname} or \cite[\textsection A]{luriestheorem}. We will not use formal geometry in any depth and one should keep in mind that classical Deligne--Mumford stacks also define formal Deligne--Mumford stacks whose topology on each \'{e}tale open is discrete.\\

Write $\M_\Ell^\sm$ for the moduli stack of smooth elliptic curves, and $\M_\Ell$ for its compactification, which has a moduli interpretation as the moduli stack of generalised elliptic curves; see \cite{delignemumford}, \cite{dr}, \cite{conrad}, or \cite{amodulardesciption} for more on such objects. Our generalised elliptic curves will \textbf{always} have irreducible geometric fibres, so either elliptic curves or N\'{e}ron 1-gons.

\begin{mydef}
The \emph{moduli stack of forms of $\G_m$} is defined as the quotient stack
\[\M_{\G_m}=(\Spec \Z)/C_2=BC_2.\]
A \emph{form of $\G_m$} over a ring $R$ is an abelian group scheme $G$ over $R$ which under a faithfully flat base change is equivalent to $\G_m$; see \cite[Pr.A.4]{tylerniko} for a proof that $\M_{\G_m}$ classifies such objects. This comes with a natural closed immersion $\M_{\G_m}\to \M_\Tate = \Spec \Z\llbracket q\rrbracket /C_2$ defined by setting $q=0$; this map plays an important role in \cite{hilllawson}, but it will not appear in this paper again.
\end{mydef}

To study the formal groups associated with generalised elliptic curves in a $p$-complete setting, we will use \emph{$p$-divisible groups} also known as \emph{Barsotti--Tate groups}. For a fixed prime $p$, write $\M_\BT$ for the \emph{moduli stack of $p$-divisible groups} and $\M_\BTn$ for the substack of $p$-divisible groups of height $n$. These are related to smooth elliptic curves and forms of $\G_m$ through the following construction of Tate \cite[\textsection2]{tatepdiv}.

\begin{mydef}\label{pdivisiblegroupsmap}
If $E$ is a smooth elliptic curve or form of $\G_m$ define the associated \emph{$p$-divisible group $E[p^\infty]$} of $E$ to have $n$th level the $p^n$-torsion subgroup $E[p^\infty]_n= E[p^n]$. This operation is functorial, and we obtain the following morphisms of stacks:
\[[p^\infty]\colon \M_\Ell^\sm\to \M_\BTtwo\qquad [p^\infty]\colon \M_{\G_m}\to \M_\BTone\]
\end{mydef}

Let $\G$ be a $p$-divisible group over a ring $R$. If $R$ is $p$-complete, there is a formal group $\G^\circ$ associated with $\G$ called its \emph{identity component}; see \cite[Th.2.0.8]{ec2name} for the construction of $(-)^\circ$ in this generality, and \cite[\textsection2.2]{tatepdiv} for the inverse functor defined for \emph{connected} $p$-divisible groups. This assignment is also compatible with the \emph{formal completion} $\widehat{(-)}$ of group schemes at their identity element, in the sense that after a base change over $\Spf\Z_p$, there is a morphism $(-)^\circ \colon\widehat{\M}_\BT\to\widehat{\M}_\FG$ of stacks and diagrams
\begin{equation}\label{twodiagramsofstacks}\begin{tikzcd}
	&	{\widehat{\M}^\sm_\Ell}\ar[rd, "{\widehat{(-)}}"]\ar[dl, "{[p^\infty]}", swap]	&	\\
{\widehat{\M}_\BTtwo}\ar[rr, "{(-)^\circ}"]	&&	{\widehat{\M}_\FG}
\end{tikzcd}\qquad
\begin{tikzcd}
	&	{\widehat{\M}_{\G_m}}\ar[rd, "{\widehat{(-)}}"]\ar[dl, "{[p^\infty]}", swap]	&	\\
{\widehat{\M}_\BTone}\ar[rr, "{(-)^\circ}"]	&&	{\widehat{\M}_\FG}
\end{tikzcd}\end{equation}
where $\M_\FG$ is the \emph{moduli stack of formal groups}; see \cite[\textsection6]{naumann07}. The commutativity of the above diagrams follows along the lines of \cite[Pr.7.4.1]{ec2name}; also see \cite[\textsection2.2.4]{ec2name}.


\subsection{Constructions using $p$-divisible groups}\label{pdivsconstuctionssection}
Fix a prime $p$, and let $C_{\BTn}$ be the subcategory of $\Fun(\CRing, \Spc)_{/\widehat{\M}_\BTn}$ spanned by those objects $\G\colon \x\to \widehat{\M}_\BTn$ where $\x$ is represented by a formal Deligne--Mumford stack of finite presentation over $\Spf \Z_p$ and $\G$ is a formally \'{e}tale morphism; see \cite[Df.2.1]{luriestheorem}. Equip this category with the \'{e}tale topology by declaring a map to be an \'{e}tale cover if the underlying map of Deligne--Mumford stacks is such. This is a particular subsite of the site $\C_{\Z_p}$ of \cite[Df.1.10]{luriestheorem}; see \cite[Pr.1.13]{luriestheorem}. The following is then a simplification of \emph{Lurie's theorem}; see \cite[Th.1.11]{luriestheorem}.

\begin{theorem}\label{luriestheorem}
Let $p$ be a prime and $n$ a positive integer. Then there is an \'{e}tale hypersheaf of $\E_\infty$-rings $\O^\top_\BTn$ on ${C}_\BTn$ such that for each affine $\G\colon \Spf R\to \widehat{\M}_\BTn$ in $C_\BTn$, the $\E_\infty$-ring $\O^\top_\BTn(\G)=\EE$ has the following properties:
\begin{enumerate}
\item $\EE$ is complex periodic
\item The groups $\pi_k \EE$ vanish for all odd integers $k$.
\item There is a chosen natural isomorphism of rings $\pi_0 \EE\simeq R$.
\item There is a chosen natural isomorphism of formal groups $\G^\circ\simeq \widehat{\G}_\EE^{\QQ_0}$ over $\Spf R$, between the identity component of $\G$ and the classical Quillen formal group of $\EE$.
\end{enumerate}
\end{theorem}

This theorem can be applied in a few concrete cases of interest to us. The following is due to Lurie, and proofs can be found in \cite[\textsection6]{luriestheorem}.

\begin{cor}\label{identificationcor}
For every prime $p$, the morphisms of stacks
\[\widehat{\M}_{\G_m}\xrightarrow{[p^\infty]} \widehat{\M}_\BTone\qquad \widehat{\M}_\Ell^\sm\xrightarrow{[p^\infty]} \widehat{\M}_\BTtwo\]
lie in $C_\BTn$, for $n=1$ and $2$, respectively. Moreover, we have an equivalence of $\E_\infty$-rings
\[\O^\top_\BTone(\widehat{\M}_{\G_m})\simeq \KO_p\]
and the diagram of $\infty$-categories
\[\begin{tikzcd}
{\left(\DM^\et_{/\M^\sm_\Ell}\right)^\op}\ar[d, "{\times \Spf \Z_p}"]\ar[rr, "{\O^\top}"]	&&	{\CAlg}\ar[d, "{(-)_p}"]	\\
{\left(\fDM^\et_{/{\widehat{\M}_{\Ell}^\sm}}\right)^\op}\ar[r, "{[p^\infty]_\ast}"]	&	{\left(C_{\BTtwo}\right)^\op}\ar[r, "{\O^\top_\BTtwo}"]	&	\CAlg
\end{tikzcd}\]
commutes, where $\O^\top$ is the Goerss--Hopkins--Miller sheaf of $\E_\infty$-rings of \cite{tmfbook} or \cite[\textsection7]{ec2name}. In particular, there is the following equivalence of $\E_\infty$-rings:
\[\O^\top_\BTtwo(\widehat{\M}_\Ell^\sm)\simeq \TMF_p\]
\end{cor}

By reformulating the above results, we obtain a functorial description of Adams operations on $\KO_p$ and $\TMF_p$. First, we need to define two sites.

\begin{mydef}\label{simplesites}
Fix a prime $p$. Define the category $C_{\sm/\BTtwo}$ as follows:
\begin{itemize}
\item Objects are étale morphisms $E\colon \xx \to \widehat{\M}_{\Ell}$ from a formal Deligne--Mumford stacks to the moduli stack of smooth elliptic curves $\M_\Ell^\sm \times \Spf \Z_p$. 
\item	Morphisms $(\xx,E)\to (\xx',E')$ given by a pair $(f,\phi)$ of a morphism of formal Deligne--Mumford stacks $f\colon \xx\to \xx'$ is a morphism of formal Deligne--Mumford stacks and $\phi\colon E[p^\infty]\simeq f^\ast E'[p^\infty]$ an isomorphism of $p$-divisible groups over $\xx$; where $[p^\infty]$ denotes the morphisms of \Cref{pdivisiblegroupsmap}.
\end{itemize}
Similarly, define a category $C_{\G_m/\BTone}$ to have an objects étale morphisms $G\colon \xx \to \widehat{\M}_{\G_m}$ defining $G$, a form of $\G_m$, and morphisms are morphisms of stacks and isomorphism of $p$-divisible groups associated with these forms of $\G_m$. Equip both of these categories with the \'{e}tale topology through the forgetful functor to formal Deligne--Mumford stacks.
\end{mydef}

These sites mirror those defined in \cite[Df.1.5]{heckeontmf}; in fact, $C_{\sm/\BTtwo} = \widehat{\mathrm{Isog}}$, using the notation of \emph{loc.\ cit}.

\begin{prop}\label{periodicitmfmotherufkcer}
Fix a prime $p$ and write $C$ for either $C_{\G_m/\BTone}$ or $C_{\sm/\BTtwo}$. There exists an \'{e}tale hypersheaf of $\E_\infty$-rings $\O^C$ on $C$ such that for an affine $E\colon\Spf R\to \widehat{\M}_{\Ell}$ in $C$, the $\E_\infty$-ring $\O^C(R)=\EE$ defines an elliptic cohomology theory for $E$, natural in $C$.
\end{prop}

Let us detail what we mean by the above elliptic cohomology theories being natural in $C$---this is simply unravelling the naturality in \Cref{luriestheorem}. Fix $C=C_{\sm/\BTtwo}$ for definiteness. For a morphism $(f,\phi)\colon (\Spf R,E)\to (\Spf R',E)$ between affine objects in $C$, then $\EE=\O^\sm(R)$ and $\EE'=\O^\sm(R')$ are natural elliptic cohomology theories, so the isomorphisms $\pi_0 \EE\simeq R$ and $\pi_0 \EE'\simeq R'$ commute with the maps $f^\ast\colon R'\to R$ and $\EE'\to \EE$. Moreover, the map $\EE'\to \EE$ induces a morphism of formal groups over $\Spf R$
\begin{equation}\label{naturalityinpdivisiblegroups}\widehat{\G}^{\QQ_0}_\EE\to f^\ast \widehat{\G}^{\QQ_0}_{\EE'}.\end{equation}
The naturality of the isomorphisms $\al\colon \widehat{E}\simeq \widehat{\G}^{\QQ_0}_\EE$ and $\al'\colon \widehat{E}'\simeq \widehat{\G}^{\QQ_0}_{\EE'}$ in $C$ means they commute with (\ref{naturalityinpdivisiblegroups}) and the morphism that $\phi\colon E[p^\infty]\to f^\ast E'[p^\infty]$ induces on formal groups by taking identity components.

\begin{proof}
There is a functor
\begin{equation}\label{pdivisiblesitfunctor}[p^\infty]\colon C\to \fDM_{/\widehat{\M}_\BTn}\end{equation}
sending a pair $(\x,E)$ to the pair $(\x,E[p^\infty])$, where $n$ depends in an obvious way on the choice of $\C$. Now we consider our two cases.
\begin{itemize}
\item We claim that the morphism of formal Deligne--Mumford stacks $\widehat{\M}_{\G_m}\to \widehat{\M}_\BTone$ is formally \'{e}tale. It would suffice to check this on the $2$-fold \'{e}tale cover $\Spf \Z_p\to \widehat{\M}_{\G_m}$, and the morphism $\Spf \Z_p\to \widehat{\M}_\BTone$ is formally \'{e}tale as it classifies the universal deformation of the multiplicative $p$-divisible group $\mu_{p^{\infty}}$ over $\F_p$; just take $\pi_0$ of \cite[Cor.3.1.19]{ec2name}. The functor $[p^\infty]$ then factors through ${C}_\BTone$. We can then define $\O^\mult$ as the following composition:
\[\O^\mult\colon C_{\G_m/\BTone}^\op\xrightarrow{[p^\infty]^\op} {C}_\BTone^\op \xrightarrow{\O^\top_\BTone} \CAlg\]
\item By the Serre--Tate theorem, see \cite{serretatevolII} for the original source and \cite[Ex.2.6]{luriestheorem} for an explanation in this context, the map $\widehat{\M}_\Ell \to \widehat{\M}_\BTtwo$ is formally étale. In particular, the above functor (\ref{pdivisiblesitfunctor}) factors through ${C}_{\BTtwo}$. This then yields a functor $[p^\infty]\colon C_{\sm/\BTtwo}\to {C}_\BTtwo$. Define $\O^\sm$ as the following composition:
\[\O^\sm\colon C_{\sm/\BTtwo}^\op\xrightarrow{[p^\infty]^\op} {C}_\BTtwo^\op \xrightarrow{\O^\top_\BTtwo} \CAlg\]
\end{itemize}
In either case, the first functor does not change the underlying formal Deligne--Mumford stack, so \'{e}tale hypercovers are sent to \'{e}tale hypercovers, and $\O^\top_\BTn$ are \'{e}tale hypersheaves, so we have two \'{e}tale hypersheaves $\O^C$ for varying $C$. These sheaves $\O^C$ satisfy the desired properties by \Cref{luriestheorem} and \Cref{identificationcor} together with the identification (\ref{twodiagramsofstacks}) of the identity component of $p$-divisible groups with the associated formal groups of smooth elliptic curves or forms of $\G_m$.
\end{proof}

In particular, we can now (re)define Adams operations on $\KO_p$ and $\TMF_p$; see \cite[\textsection6.4]{luriestheorem} for a previous formulation, more properties, and the relation to classical Adams operations.

\begin{mydef}\label{adamsoperationsofperiodics}
For each prime $p$ and each $p$-adic unit $k\in\Z_p^\times$, define the (auto)morphisms of $\E_\infty$-rings
\[\psi^k\colon \KO_p\to \KO_p\qquad \psi^k\colon \TMF_p\to \TMF_p\]
by applying $\O^C$ of \Cref{periodicitmfmotherufkcer} to the $k$-fold multiplication map of $p$-divisible groups associated with the universal group schemes over $\widehat{\M}_{\G_m}$ and $\widehat{\M}_\Ell^\sm$, respectively---note this $k$-fold multiplication is an equivalence of $p$-divisible groups as $k\in \Z_p^\times$.
\end{mydef}

\begin{remark}
From knowledge about $\KU_p$ as a Lubin--Tate theory of height $1$ at the prime $p$, we know that the above Adams operations on $\KO_p$ are all such automorphisms of this $\E_\infty$-ring; see \cite[\textsection7]{gh04} or \cite[\textsection5]{ec2name}. There are no other obvious $\E_\infty$-automorphisms of $\TMF_p$, at least to the author. After $K(2)$-localisation, there is also much of a height $2$-Morava stabiliser group action on $\TMF_p$, and perhaps some of these automorphisms can be lifted to the $E_2$-local $\TMF_p$. The reader interested in initiating such lifts should start with $p=3$ and the $K(2)$-local discussions of $\TMF_p$ found in \cite{ktwospheremark,ktwolocalsphere}.
\end{remark}

From these operations, we also obtain Adams operations on \emph{Tate $K$-theory}; an exploration of Tate $K$-theory through the lens of spectral algebraic geometry can be found in \cite{globaltate}.

\begin{mydef}\label{adamsoperationsontate}
For an $\E_\infty$-ring $A$, we define $A\llbracket q\rrbracket$ as the completion of $A \otimes \Sigma_+^\infty \mathbf{N}$ at the element $q\in \pi_0 A \otimes \Sigma_+^\infty \mathbf{N} \simeq \pi_0 A[q]$; the isomorphism here comes from a degenerating Tor-SS from the flatness of $\Sigma_+^\infty \mathbf{N}$ over $\Sph$. Notice that the natural map of $\E_\infty$-rings
\[\KO\llbracket q \rrbracket_p \xrightarrow{\simeq} \KO_p\llbracket q \rrbracket_p\]
is an equivalence. Indeed, using standard facts about $p$-completion and computing the homotopy groups of each side, this boils down to the classical fact that $\Z\llbracket q\rrbracket \to \Z_p\llbracket q\rrbracket$ induces an isomorphism on classical $p$-completions; a fact that is obvious as we have natural identifications
\begin{equation}\label{littleadiccomputation}\Z\llbracket q\rrbracket/p^n \Z\llbracket q\rrbracket \simeq (\Z/p^n)\llbracket q\rrbracket \simeq \Z_p\llbracket q\rrbracket/p^n \Z_p\llbracket q\rrbracket.\end{equation}
In particular, $\KO\llbracket q\rrbracket_p$ comes equipped with Adams operations $\psi^k$ for each $k\in \Z_p^\times$ from those on $\KO_p$ from \Cref{adamsoperationsofperiodics}. Similarly, the natural map of $\E_\infty$-rings
\[\KO\llparenthesis q\rrparenthesis_p \xrightarrow{\simeq} \KO_p\llparenthesis q\rrparenthesis_p\]
is also an equivalence, as the quotients of (\ref{littleadiccomputation}) commute with inverting $q$, a type of colimit. In particular, the $\E_\infty$-ring $\KO\llparenthesis q\rrparenthesis_p$ can be equipped with Adams operations $\psi^k$ for each $k\in \Z_p^\times$ from those on $\KO_p$.
\end{mydef}

It is clear from the above definitions that all of the maps of $\E_\infty$-rings
\[\KO_p\to \KO\llbracket q\rrbracket_p \to \KO_p,\qquad \KO_p \to \KO\llparenthesis q\rrparenthesis_p \to \KO\llbracket q\rrbracket_p,\]
where the second coming from setting $q=0$, all naturally commute with the each Adams operation $\psi^k$.


\subsection{Proof of \Cref{thaversiontwo}}\label{gluingtogether}

To glue together our Adams operations on $\TMF_p$ with those on $\KO\llbracket q\rrbracket_p$ will use Goerss--Hopkins obstruction theory. This will destroy much of our functorality, only allowing us to construct each Adams operation $\psi^k\colon \Tmf_p\to \Tmf_p$ for each $k\in \Z_p^\times$ in isolation. We encourage the reader to remind themselves of the $K(1)$-local stable homotopy theory of \cite[\textsection 12.6-7]{tmfbook}. We would again like to thank an anonymous referee for suggestions to simplify this section.\\

For this subsection, fix a prime $p$ and a $p$-adic unit $k\in \Z_p^\times$.

\begin{prop}\label{thalemmaone}
There is a morphism of $\E_\infty$-rings $\Lambda\colon \TMF_p\to \KO\llparenthesis q\rrparenthesis_p$ such that the diagram of $\E_\infty$-rings
\begin{equation}\label{diagramwithsomelambs}
\begin{tikzcd}
{\TMF_p}\ar[d, "{\Lambda}"]\ar[r,"{\psi^k}"]	&	{\TMF_p}\ar[d, "{\Lambda}"]	\\
{\KO\llparenthesis q\rrparenthesis_p}\ar[r, "{\psi^k}"]			&	{\KO\llparenthesis q\rrparenthesis_p}
\end{tikzcd}\end{equation}
 commutes up to homotopy, where the horizontal maps are the Adams operations from \Cref{adamsoperationsofperiodics,adamsoperationsontate}.
\end{prop}

Of course, this morphism $\Lambda$ is (up to homotopy) the ``evaluation at the cusp'' map $\Tmf\to \KO\llbracket q\rrbracket$ with $\Delta^{24}$ inverted and $p$-completed; see \cite[\textsection A]{hilllawson} for an obstruction theoretic approach to this map and \cite{globaltate} for a spectral algebro-geometric approach. This follows from the uniqueness of this map up to $1$-homotopy, a consequence of the methods of \cite[Pr.A.6]{hilllawson} or the proof below.


\begin{proof}
First, note that $\KO\llparenthesis q\rrparenthesis_p$ is $K(1)$-local, as it is a $p$-complete $\KO$-module (see \cite[Rmk.A.2]{hilllawson}), so we may $K(1)$-localise $\TMF_p$ and work in $\CAlg_{K(1)}$. We will now use $K(1)$-local Goerss--Hopkins obstruction theory, as found in \cite[\textsection12.7]{tmfbook}, \cite[\textsection A]{hilllawson}, or \cite[\textsection5.4]{lawsonnaumannbp}, for example. The following arguments depend on the parity of $p$.

\paragraph{(For $p\neq 2$)} Recall the $p$-adic $K$-theory $\K^\wedge_\ast R$ of an $\E_\infty$-ring $R$ is defined as $\pi_\ast L_{K(1)}(\KU\otimes R)$, and comes equipped with the structure of a $\theta$-algebra; see \cite[\textsection12.6]{tmfbook} or \cite{gh04}. It follows from the arguments of \cite[Pr.A.6]{hilllawson} that the $p$-adic $K$-theory functor induces the following bijection of sets:
\begin{equation}\label{mappingcomponents}\pi_0 \Map_{\CAlg_{K(1)}}\left(L_{K(1)}\TMF, \KO\llparenthesis q\rrparenthesis_p\right)\xrightarrow{\simeq} \Hom_{\theta\Alg_{(\KU_p)_\ast}}\left(\K^\wedge_\ast \TMF,\K^\wedge_\ast \KO\llparenthesis q\rrparenthesis\right)\end{equation}
Indeed, this is due to the isomorphism of $\theta$-algebras $\K^\wedge_\ast\TMF\simeq (\KU_p)_\ast\otimes_{\Z_p} V$, where $V$ is the $p$-adic ring representing \emph{smooth} elliptic curves $E$ with a chosen isomorphism between $\widehat{E}$ and $\widehat{\G}_m$ (see \cite[\textsection12.5]{tmfbook} for a discussion of $V$, which is the smooth variant of what is written there as $V^\wedge_\infty$), and the fact that this $V$ is formally smooth over $\Z_p$; see \cite[Lm.12.7.9]{tmfbook}. By (\ref{mappingcomponents}), we see it suffices to study the $p$-adic $K$-theory of $\TMF$ and $\KO\llparenthesis q\rrparenthesis$. Following \cite[Pr.A.4]{hilllawson}, we can also calculate the $\theta$-algebra $\K^\wedge_\ast \KO\llparenthesis q\rrparenthesis$ as $(\KU_p)_\ast\otimes_{\Z_p}V_\Tate$, where $V_\Tate$ is now defined as the universal $p$-adic $\Z\llparenthesis q\rrparenthesis$-algebra with an isomorphism class of pairs of an invariant $1$-form on the smooth Tate curve $T$ and a chosen isomorphism between $\widehat{T}$ and $\widehat{\G}_m$.\footnote{Note that this calculation of $\K^\wedge_\ast\KO\llparenthesis q\rrparenthesis$ holds for all primes, as the arguments calculating $\K^\wedge_\ast\KO\llbracket q\rrbracket$ from \cite[Pr.A.4]{hilllawson} also hold in this generality.} There is a canonical map $\lambda\colon V\to V_\Tate$ of $p$-adic rings as the smooth Tate curve $E$ is an elliptic curve, and as explained in the proof of \cite[Pr.4.49]{hilllawson}, this morphism defines a map of $\theta$-algebras. The map $\lambda$ is known as the $p$-adic \emph{$q$-expansion map}, and by (\ref{mappingcomponents}), we can recognise this map by our desired morphism of $K(1)$-local $\E_\infty$-rings $\Lambda \colon L_{K(1)}\TMF\to \KO\llparenthesis q\rrparenthesis_p$. We are now required to show the diagram of $\E_\infty$-rings (\ref{diagramwithsomelambs}) commutes. Appealing to (\ref{mappingcomponents}) again, we are reduced to show that the above diagram commutes after applying $p$-adic $K$-theory. As we know the $p$-adic $K$-theory of all of the above $\E_\infty$-rings, and these $p$-adic $K$-theories are all base changed from their zeroth $p$-adic $K$-theory, it suffices to check (\ref{diagramwithsomelambs}) commutes after applying zeroth $p$-adic $K$-theory. We will get back to this shortly, once we bring the case for even $p$ up to speed.

\paragraph{(For $p=2$)} Recall from \cite[Df.12.7.10]{tmfbook} (as well as the appendix of that chapter) that the $2$-adic real $K$-theory $\KO^\wedge_\ast R$ of an $\E_\infty$-ring $R$ is defined as $\pi_\ast L_{K(1)}(\KO\otimes R)$, and naturally has the structure of a \emph{reduced} graded $\theta$-algebra, meaning that $\psi^{-1}$ acts trivially. There is a form for Goerss--Hopkins obstruction theory in this situation for \emph{Bott periodic} $\E_\infty$-rings. An $\E_\infty$-ring $R$ is said to be \emph{Bott periodic} if $\K^\wedge_\ast R$ is torsion-free and concentrated in even degrees, and the natural map $\KO^\wedge_0 R\to \K^\wedge_0 R$ is an isomorphism. Bott periodic Goerss--Hopkins obstruction theory then states that if $R_1$ and $R_2$ are two $K(1)$-local Bott periodic $\E_\infty$-rings and $f_\ast\colon \KO^\wedge_\ast R_1\to \KO^\wedge_\ast R_2$ is a morphism of reduced graded $\theta$-algebras, then the obstructions to the lifting $f_\ast$ to a map $f\colon R_1\to R_2$ of $K(1)$-local $\E_\infty$-rings lie in the following Andr\'{e}--Quillen cohomology groups:
\begin{equation}\label{firstobstructiongroups}H^s_{\theta\Alg_{(\KO_2)_\ast}^\red}(\KO^\wedge_\ast R_1,\KO^\wedge_\ast R_2[-s+1])\qquad s\geq 2\end{equation}
Moreover, obstructions to the uniqueness of $f$ recognising $f_\ast$ up to homotopy live in the following cohomology groups:
\begin{equation}\label{secondobstructiongroups}H^s_{\theta\Alg_{(\KO_2)_\ast}^\red}(\KO^\wedge_\ast R_1,\KO^\wedge_\ast R_2[-s])\qquad s\geq 1\end{equation}
We claim that for $R_1=L_{K(1)}\TMF$ and $R_2=\KO\llparenthesis q\rrparenthesis_2$, both families of obstruction groups above vanish. To show this, consider the vanishing criteria of \cite[Lm.12.7.13]{tmfbook}:
\begin{enumerate}
\item[(1)]	The $\E_\infty$-ring $\tmf$ is Bott periodic by construction (\cite[Rmk.12.7.12]{tmfbook}) and it follows that its localisation $\TMF$ is also Bott periodic. To see $\KO\llparenthesis q\rrparenthesis$ is Bott periodic, we first refer to the calculation that $\K^\wedge_\ast \KO\llparenthesis q\rrparenthesis$ is isomorphic to the $\theta$-algebra $(\KU_p)_\ast\otimes_{\Z_p} V_\Tate$ discussed above under the assumption that $p$ is odd---indeed, this calculation holds for all primes. It is rather formal that the natural map
\begin{equation}\label{secondbpcondition}\KO_0^\wedge\KO\llparenthesis q\rrparenthesis\to \K^\wedge_0 \KO\llparenthesis q\rrparenthesis\end{equation}
is an isomorphism. Indeed, as the map of $\E_\infty$-rings $\KO\to \KO\llparenthesis q\rrparenthesis$ is flat, then for any $\KO$-module $M$ we obtain natural isomorphisms
\[M_\ast \KO\llparenthesis q\rrparenthesis=\pi_\ast(M\otimes \KO\llparenthesis q\rrparenthesis)\simeq \pi_\ast(M\otimes \KO\otimes_{\KO}\KO\llparenthesis q\rrparenthesis)\simeq \pi_\ast(M\otimes \KO)\otimes \Z\llparenthesis q\rrparenthesis\]
by a degenerating K\"{u}nneth spectral sequence. For $M=\KO$ or $\KU$, we then see that the map (\ref{secondbpcondition}) is the base change of the classical isomorphism $\KO_0\KO\simeq \K_0\KO$ over $\Z\llparenthesis q\rrparenthesis$, and then $2$-completed, and hence is an isomorphism.
\item[(2)]	The mod $2$-reduction of $V^{\Z_2^\times}$ is formally smooth over $\F_2$. Indeed, the $\Z_2^\times$-action factors through a $\Z_2^\times/\{\pm1\}$-action as $[-1]$ acts trivially on $V$; see \cite[Lm.12.7.14(1)]{tmfbook}. The $\Z_2^\times/\{\pm1\}$-fixed points of $V$ are $V_2$, the $2$-adic ring representing the moduli stack $\M_\Ell^\ord(4)$ whose $S$-points consist of a pair of a smooth elliptic curve $E$, with ordinary mod $2$-reduction, and level structure given by an isomorphism of finite group schemes $\mu_{4}\simeq \widehat{E}[4]$; this is \cite[Lm.12.7.14]{tmfbook}. The (affine) stack $\M_\Ell^\ord(4)\otimes \F_2$ is smooth over $\F_2$, which proves our claim.
\item[(3)]	To see the continuous cohomology groups $H^s_c(\Z_2^\times/\{\pm 1\},V_\Tate/2 V_\Tate)$ vanish for $s\geq 1$, it suffices to see that $\Spf V_\Tate$ is an ind-Galois torsor for the group $\Z_2^\times/\{\pm 1\}$ over $\Spf \Z\llparenthesis q\rrparenthesis_2$. This follows by observing that $V_\Tate$ can be explicitly written (\`{a} la \cite[Pr.A.4]{hilllawson}) as the set of continuous maps from $\Z_2^\times/\{\pm 1\}$ into $\Z\llparenthesis q\rrparenthesis_2$, where $\Z_2^\times/\{\pm1\}$ acts by conjugation, and that using this expression for $V_\Tate$ its $\Z_2^\times/\{\pm1\}$-fixed points are precisely $\Z\llparenthesis q\rrparenthesis_2$.
\item[(4)]	In part (2) above, we saw $V^{\Z_2^\times/\{\pm1\}}$ is given by the $2$-adic ring $V_2$, so $V_2\to V$ is a $\Z_2^\times/\{\pm1\}$-ind-Galois extension. In particular, $V_2\to V$ is ind-\'{e}tale. By base change, we see that the mod $2$-reduction of this inclusion of fixed points is also ind-\'{e}tale.
\end{enumerate}

The four conditions above line up with the four hypotheses of \cite[Lm.12.7.13]{tmfbook}, and we consequently see that the obstruction groups (\ref{firstobstructiongroups}-\ref{secondobstructiongroups}) all vanish. Hence the commutativity of (\ref{diagramwithsomelambs}) can be checked on $2$-adic $\KO$-homology. By \cite[Lm.12.7.11]{tmfbook}, we see the $2$-complete $\KO$-theory of a Bott periodic $\E_\infty$-ring $R$ naturally depends on its zeroth $2$-complete $K$-theory:
\[\KO^\wedge_\ast R\simeq (\KO_2)_\ast\otimes_{\Z_2}\K^\wedge_0 R\]
Consequently, just like in the case for an odd prime $p$, we are reduced to studying the zeroth $p$-adic $K$-theory of (\ref{diagramwithsomelambs}).

\paragraph{(Back to general $p$)} It suffices to show that (\ref{diagramwithsomelambs}) commutes after applying zeroth $p$-adic $K$-theory:
\[\begin{tikzcd}
{V}\ar[r, "{\psi^k_\sm}"]\ar[d, "{\lambda}"]	&	{V}\ar[d, "{\lambda}"]	\\
{V_\Tate}\ar[r, "{\psi^k_\Tate}"]		&	{V_\Tate}
\end{tikzcd}\]
By construction, the morphism $\lambda$ is one of $\theta$-algebras with respect to the \emph{algebraic Adams operations} on both $V$ and $V_\Tate$, given by on $S$-valued points of $\Spf V$ by $(E,\al)\mapsto (E,\al\circ [k])$. Hence it suffices to show that the $p$-adic $K$-theory of the operations $\psi^k_\sm$ and $\psi^k_\Tate$ agree with the relevant algebraic Adams operations. This will follow from our construction of these operations from \textsection\ref{pdivsconstuctionssection}.\\

Let us begin with the $\psi^k_\sm$-case. For any object $X$ in $C_{\sm/\BTtwo}$, we can define $\psi^k_X\colon \O^\sm(X)\to \O^\sm(X)$ by applying $\O^\sm$ to the $k$-fold multiplication maps on the associated $p$-divisible groups; see \Cref{periodicitmfmotherufkcer} and \Cref{adamsoperationsofperiodics}, or \cite[Df.6.16]{luriestheorem}. If $X=\Spf R$ is affine, then \cite[Lm.12.6.1]{tmfbook} supplies us with the morphisms and equivalences of formal stacks
\begin{equation}\label{pullbackfrommark}\begin{tikzcd}
{\Spf \K^\wedge_0 A=\Spf V_{R}}\ar[r]\ar[d]	&	{\M_\Ell^{\sm,\ord}(p^\infty)=\Spf V=\Spf \K^\wedge_0 \TMF}\ar[d]\ar[r]	&	{\M_\Ell^\ord(p^\infty)=\Spf \K^\wedge_0 \tmf}\ar[d]	\\
{\Spf R^\ord}\ar[r]					&	{\M_\Ell^{\sm,\ord}}\ar[r]	&	{\M_\Ell^\ord}
\end{tikzcd}\end{equation}
where $A=\O^\sm(X)$, $\Spf R^\ord$ is $\Spf R$ base-changed over $\M_\Ell^\ord\to \widehat{\M}_\Ell$, and the square above is Cartesian. The naturality with respect to $C_{\sm/\BTtwo}$ of the isomorphism in condition (iv) of \Cref{luriestheorem} shows that the map that $\psi^k_X$ induces the $k$-fold multiplication map on the associated Quillen formal groups. Hence the map $\psi^k_X\colon V_{R}\to V_{R}$ is represented by the pair $(E,[k]\circ \al)$, where $E$ is the universal smooth elliptic curve over $V$ pulled back to $V_{R}$, $\al$ is the base change of the universal isomorphism $\widehat{E}\simeq\widehat{\G}_m$ to $V_{R_\sm}$, and $[k]$ is the $k$-fold multiplication map on formal groups---such a map of $p$-divisible groups induced such a map on formal groups.\\

When $X=\widehat{\M}_\Ell^\sm$, then we can choose an affine \'{e}tale cover $\Spf R\to \widehat{\M}_\Ell^\sm$ and again consider the diagram of formal stacks (\ref{pullbackfrommark}). In this case, the lower-horizontal map is faithfully flat by assumption, so the upper-horizontal map is also faithfully flat. In particular, the map of rings $V\to V_R$ is injective. From the argument above, the algebraic Adams operations on $V_R$ and those induced by $\psi^k_R$ agree. Moreover, the map $V\to V_R$ is induced by $\tmf_p\to \O^\sm(\Spf R)$, hence it commutes with the Adams operations induced by $\psi^k_\sm$ and $\psi^k_R$. Finally, the map $V\to V_R$ also commutes with the algebraic Adams operations as we again appeal to \cite[Pr.4.49]{hilllawson} which states that this holds if $V\to V_R$ is a map of rings over $\widehat{\M}_\Ell$. From these facts and the injectivity of $V\to V_R$, we see that the algebraic Adams operations on $V$ agree with those induced by $\psi^k_\sm$.\\

The $\psi^k_\Tate$-case is analogous. Indeed, coping the above affine argument for $\KU_p$, we see that the Adams operations $\psi^k$ on $\KU_p$ induce the algebraic Adams operations on $p$-adic $K$-theory. As the map of $\E_\infty$-rings $\KO\to \KU$ induces an isomorphism of zeroth $p$-adic $K$-theory, we see that the operations $\psi^k$ on $\KO_p$, themselves induced from $\psi^k$ on $\KU_p$, also induce the algebraic operations on zeroth $p$-adic $K$-theory. As the Adams operations on $\KO\llparenthesis q\rrparenthesis_p$ are determined by those on $\KO_p$, and that likewise the algebraic Adams operations on the zeroth $p$-adic $K$-theory of $\KO\llparenthesis q\rrparenthesis_p$ are determined by those on $\KO_p$, we obtain the desired result.
\end{proof}

For primes $p\neq 2$, there is a strengthening of the previous proposition. Let us equip $\TMF_p$, $\KO_p$, $\KO\llbracket q\rrbracket_p$, and $\KO\llparenthesis q\rrparenthesis_p$ with an $\F_p^\times$-action using the coherent Adams operations of \Cref{adamsoperationsofperiodics,adamsoperationsontate} together with the multiplicative lift $\F_p^\times\leq \Z_p^\times$. From these definitions, the Adams operations $\psi^k$ are $\F_p^\times$-equivariant as automorphisms of $\E_\infty$-rings.

\begin{prop}\label{thalemmathree}
There is a morphism of $\E_\infty$-rings $\Lambda\colon \TMF_p\to \KO\llparenthesis q\rrparenthesis_p$ with $\F_p^\times$-action such that the diagram (\ref{diagramwithsomelambs}) inside $\CAlg^{B\F_p^\times}$ commutes up to homotopy.
\end{prop}

\begin{proof}
One can carry out the whole argument used to prove \Cref{thalemmaone} for odd primes in the setting of $\F_p^\times$-equivariant $\E_\infty$-rings, so the category $\CAlg^{B \F_p^\times}$. As discussed in \cite[\textsection5.1]{vesnaduality}, there is a $G$-equivariant form of Goerss--Hopkins obstruction theory for a finite group $G$. If the order of our group $\F_p^\times$ is not divisible by $p$, then all of the $\F_p^\times$-equivariant obstruction groups can be calculated as the $\F_p^\times$-fixed points of the non-equivariant obstruction groups used in the proof of \Cref{thalemmaone}. This allows us to run all of the arguments of \Cref{thalemmaone} in the $\F_p^\times$-equivariant setting and obtain our desired result. One can alternatively construct these $\F_p^\times$-equivariant Adams operations using the variant of Goerss--Hopkins obstruction theory found in \cite{realspectra}.
\end{proof}

We can now construct the $p$-adic Adams operations on $\Tmf_p$.

\begin{proof}[Proof of \Cref{thaversiontwo}]
Define an $\E_\infty$-ring $\Tmf_p$ using \Cref{thalemmaone} (or an $\F_p^\times$-equivariant $\E_\infty$-ring using \Cref{thalemmathree} at odd primes) via the Cartesian diagram
\[\begin{tikzcd}
{\Tmf_p}\ar[r]\ar[d]		&	{\KO\llbracket q\rrbracket_p}\ar[d]	\\
{\TMF_p}\ar[r, "{\Lambda}"]	&	{\KO\llparenthesis q\rrparenthesis_p;}
\end{tikzcd}\]
this diagram is one definition of $\Tmf_p$; see \cite{globaltate}. By (the proof of) \Cref{thalemmaone} we see that the map $\Lambda$ agrees with the usual smooth $q$-expansion map up to $1$-homotopy (as they both have the same effect on zeroth $p$-adic $K$-theory by construction), hence the pullback is homotopy equivalent to any other $\E_\infty$-ring one might call $\Tmf_p$. Moreover, \Cref{thalemmaone} equips $\Tmf_p$ with an endomorphism of $\E_\infty$-rings $\psi^k$ which agrees with the action of $\psi^k$ when restricted to $\TMF_p$ and $\KO\llbracket q\rrbracket_p$. Moreover, when $p$ is odd, \Cref{thalemmathree} constructs $\Tmf_p$ as an $\E_\infty$-ring with $\F_p^\times$-action equipped with an $\F_p^\times$-equivariant morphism of $\E_\infty$-rings $\psi^k$.
\end{proof}

Let us reiterate: these Adams operations $\psi^k$ on $\Tmf_p$ have no obvious compatibility as $k$ varies---when working over $\M_\Ell^\sm$, we have natural homotopies $\psi^k\psi^\ell\simeq \psi^{k\ell}$, for example; see \cite[Pr.6.17]{luriestheorem} or \cite[Ths.C-D]{heckeontmf}. For odd primes $p$, one can show there are homotopies between the Adams operations on $\Tmf_p$ of the form $\psi^k\psi^\ell\simeq \psi^{k\ell}$, and such homotopies are associative up to $3$-homotopy; see \cite[Th.C \& Th.3.16]{realspectra}.

\begin{remark}
One might hope that other constructions on $\TMF_p$ made possible using $p$-divisible groups also have analogues for $\Tmf_p$. For example, the morphisms $q^\ast\colon \TMF_p\to \TMF_0(\ell)_p$ defined for a prime $\ell$ distinct from $p$. The construction of these morphisms over $\M_\Ell^\sm$ is simple, they send a pair $(E,H)$ to the quotient $E/H$, but over the compactification require a lot of care; see \cite[\textsection4.4.3]{conrad} and \cite[\textsection4.7]{amodulardesciption}. Following the recipe above, one can construct morphisms $q^\ast\colon \Tmf_p\to \Tmf_0(\ell)_p$ which restrict to the above morphisms of periodic topological modular forms. There are at least two reasons one might like such additional morphisms surrounding $\Tmf$: to construct Hecke operators on $\Tmf_p$ and hence also $\tmf_p$, akin to those on $\TMF_p$ found in \cite[\textsection2]{heckeontmf}, and to construct connective versions of Behrens' $Q(N)$ spectra of \cite{ktwospheremark}. Both of these constructions will appear in future work.
\end{remark}


\subsection{Proof of \Cref{tha}}\label{newintegralstuff}

To construct the desired map of $\E_\infty$-rings $\psi^k\colon\Tmf[\frac{1}{k}]\to \Tmf[\frac{1}{k}]$ for a positive integer $k$, we will glue the morphisms $\psi^k$ on $\Tmf_p$ for each primes $p$ not dividing $k$ together with a morphism $\psi^k$ on $\Tmf_\Q$ that we will construct shortly---the techniques used here are standard.\\

Recall that $\tmf_\Q$, the rationalisation of $\tmf$, is \emph{formal} as a rational cdga;\footnote{Here, we are implicitly using the symmetric monoidal Schwede--Shipley equivalence of $\infty$-categories $\Mod_{\Q}\simeq \mathcal{D}(\Q)$; see \cite{schwedeshipleystablemodulecategoresiarecategirofmodules} or \cite[Th.7.1.2.13]{haname}.} see \cite[Pr.4.47]{hilllawson}, for example. This means that $\tmf_\Q$ is equivalent to the connective formal rational cdga $A_\ast=\Lambda_\Q[c_4,c_6]$ defined by the free $\E_\infty$-$\Q$-algebra on elements $c_4\in A_{8}$ and $c_6\in A_{12}$. Write $\Delta$ for the element $\frac{c_4^3-c_6^2}{1728}\in A_{24}$. Consider the following Cartesian square of rational cdgas:

\begin{equation}\label{rationalsquare}\begin{tikzcd}
{\Tmf_\Q}\ar[r]\ar[d]				&	{\tmf_\Q[c_4^{-1}]}\ar[d]	\\
{\tmf_\Q[\Delta^{-1}]=\TMF_\Q}\ar[r]		&	{\tmf_\Q[c_4^{-1},\Delta^{-1}]}
\end{tikzcd}\end{equation}

Fix an integer $k$. Define the endomorphism of rational cdgas $\psi^k\colon \tmf_\Q\to \tmf_\Q$ by sending $c_4$ to $k^4c_4$ and $c_6$ to $k^6c_6$. As this induces compatible endomorphisms on all the cdgas in (\ref{rationalsquare}) we obtain an endomorphism of rational cdgas $\psi^k$ on $\Tmf_\Q$. We are now ready to glue this endomorphism on $\Tmf_\Q$ to those of \Cref{thaversiontwo}.\\

We will use \Cref{thb} to prove \Cref{tha}, however, the statement of \Cref{thb} only involves the operations from \Cref{thaversiontwo}. Moreover, we only use the weaker rational version of \Cref{thb}, which can also be proven using the simpler techniques of \Cref{adamsopessections}.

\begin{proof}[Proof of \Cref{tha}]
Fix an integer $k$. Write $X$ for any $\E_\infty$-ring in the set
\[\{\TMF[\frac{1}{k}],\tmf[c_4^{-1},\frac{1}{k}],\TMF[c_4^{-1},\frac{1}{k}]\}.\]
Each each such $X$, there is the following Cartesian arithmetic fracture square of $\E_\infty$-rings:
\begin{equation}\label{arithmeticfracturesquare}\begin{tikzcd}
{X}\ar[r]\ar[d]	&	{\prod_{p\nmid k}X_p}\ar[d]	\\
{X_\Q}\ar[r,"{\al}"]	&	{\left(\prod_{p\nmid k}X_p\right)_\Q}
\end{tikzcd}\end{equation}
By \Cref{thaversiontwo} and our discussion of $\Tmf_\Q$ above, $X_p$ and $X_\Q$ both have an endomorphism $\psi^k$. Moreover, the right vertical morphism of (\ref{arithmeticfracturesquare}) commutes with these Adams operations, so to obtain Adams operations on $X$, we only have to show that the lower horizontal map commutes with Adams operations. This is easy though, as the space of $\E_\infty$-morphisms out of $X_\Q$ into a rational $\E_\infty$-ring $R$ is (a component of) the space
\[\Omega^{\infty+8}R\times \Omega^{\infty+12}R\]
as $\tmf_\Q$ is a free rational $\E_\infty$-ring and $X$ is a localisation. In particular, we see that two morphisms $\al\circ\psi^k$ and $\psi^k\circ \al$ agree up to homotopy if their images of $c_4$ and $c_6$ agree in the homotopy groups of $(\prod_{p\nmid k}X_p)_\Q$. \Cref{thb} allows us to compare our $p$-complete calculations of $\psi^k$ to the rational calculations (which follow by definitions), and we see that $\al\circ\psi^k$ and $\psi^k\circ \al$ do agree on homotopy groups, so we obtain endomorphisms of $\E_\infty$-rings $\psi^k\colon X\to X$. To glue together endomorphisms, consider the diagram of rational $\E_\infty$-rings
\begin{equation}\label{finalrationalgluing}\begin{tikzcd}
{\T}\ar[rr]\ar[dr]\ar[dd]	&&{\T[c_4^{-1}]}\ar[dd]\ar[rd]	&&	{\t[c_4^{-1}]}\ar[dd]\ar[rd]\ar[ll]	&	\\
	&	{\T}\ar[rr, crossing over]&&{\T[c_4^{-1}]}		&&	{\t[c_4^{-1}]}\ar[dd]\ar[ll, crossing over]	\\
{\prod\T_p}\ar[rr]\ar[dr]		&&{\prod\T[c_4^{-1}]_p}\ar[rd]		&&	{\prod\t[c_4^{-1}]_p}\ar[rd]\ar[ll]		&	\\
	&	{\prod\T_p}\ar[rr]\ar[from=uu, crossing over]	&&{\prod\T[c_4^{-1}]_p}\ar[from=uu, crossing over]			&&	{\prod\t[c_4^{-1}]_p}\ar[ll]
\end{tikzcd}\end{equation}
where $\T=\TMF$ and $\t=\tmf$, all diagonal maps are the respective Adams operations $\psi^k$, the products are taken over all primes $p$ not dividing $k$, and we have suppressed rationalisation everywhere. Repeating our arguments above, we see that each face in the above diagram commutes, up to a homotopy. We can then use the 2-skeleton of the left cube above to construct a map of spaces
\begin{equation}\label{finalobstruction}S^1\to \Map_{\CAlg_\Q}(\T,\prod \T[c_4^{-1}]_p)\subseteq \Omega^{\infty+8}\prod \T[c_4^{-1}]_p\times \Omega^{\infty+12}\prod \T[c_4^{-1}]_p\end{equation}
which encodes how these six homotopies (each represented above by whiskering a face in the left cube of (\ref{finalrationalgluing})) relate the six compositions from $\T$ to $\prod \T[c_4^{-1}]_p$ from (\ref{finalrationalgluing}). Note that the second map in (\ref{finalobstruction}) is the inclusion of a component. From (\ref{finalobstruction}), we see the obstruction to lifting the $2$-skeleton of the left cube of (\ref{finalrationalgluing}) in $\CAlg_\Q$ to the whole cube lies in $\pi_1$ of the codomain of (\ref{finalobstruction}), based at any choice of map $\T\to \prod\T[c_4^{-1}]_p$ displayed in (\ref{finalrationalgluing}). We see $\pi_i$ of the codomain of (\ref{finalobstruction}) vanishes for $i=1,2,3$, hence we see that the left cube of (\ref{finalrationalgluing}) admits a lift to a diagram in $\CAlg_\Q$. The same argument applies to the right cube of (\ref{finalrationalgluing}) \emph{mutatis mutandis}. Taking pullbacks along the horizontal cospans in (\ref{finalrationalgluing}) gives us the left square in the commutative diagram of $\E_\infty$-rings
\[\begin{tikzcd}
{\Tmf_\Q}\ar[d,"{\psi^k}"]\ar[r]		&	{\left(\prod_{p\nmid k}\Tmf_p\right)_\Q}\ar[d,"{\psi^k}"]	&	{\prod_{p\nmid k}\Tmf_p}\ar[d, "{\psi^k}"]\ar[l]	\\
{\Tmf_\Q}\ar[r]				&	{\left(\prod_{p\nmid k}\Tmf_p\right)_\Q}			&	{\prod_{p\nmid k}\Tmf_p}\ar[l]
\end{tikzcd}\]
and the right square commutes by definition. Taking pullbacks along the horizontal cospans again yields a morphism of $\E_\infty$-rings $\psi^k\colon \Tmf[\frac{1}{k}]\to \Tmf[\frac{1}{k}]$ whose $p$-completion at any $p$ not dividing $k$ is the Adams operation of \Cref{thaversiontwo}.
\end{proof}


\section{Calculations}\label{adamsopessections}

One can now take \Cref{tha,thaversiontwo} for granted, ie, the existence of Adams operations on $\Tmf$, as in this section, we prove \Cref{thb} (repeated below) using different techniques.

\begin{theorem}[{\Cref{thb}}]
For every odd prime $p$, every $p$-adic unit $k\in \Z_p^\times$, and every $x\in \pi_d\Tmf_p$ with $d$ positive, we have the equality
\[\psi^k(x)=\begin{cases}	x & x\in\Tors_d \\ k^{\ceil{\frac{d}{2}}}x & x\in\Free_d \end{cases}\]
where $\Tors_d\subseteq\pi_d\Tmf_p$ is the subgroup of torsion elements and $\Free_d$ is the orthogonal subgroup of \Cref{elementsintmf}. At the prime $p=2$, the above equalities hold for all $d$, however, they are vacuous for positive $d$ which are congruent to $60$ or $156$ modulo $192$.
\end{theorem}

In \Cref{achoiceofgenerators}, we define $\Free$ and relate them to the notation of Bruner--Rognes \cite{brunerrognes}. In \Cref{andersondualitysection}, a homotopical self-duality for $\Tmf$ is discussed, originally proven by Stojanoska. In \Cref{generaloperations}, we prove \Cref{thb} and on the way gather evidence for a conjecture concerning dual endomorphisms of self-dual spectra. \\

Throughout this section, we will freely use the notation of \cite[\textsection13]{tmfbook} and \cite{bauer} to indicate elements in $\pi_\ast\tmf$, and \cite{konter} for elements in $\pi_\ast\Tmf$. There are more details for the descent spectral sequence for $\Tmf$ given in \cite{smfcomputation}, where this topic is treated with synthetic spectra. Although the pictures in \cite[\textsection13]{tmfbook} are arguably the most readable, they can be misleading, for example, the vertical axis is neither the Adams or Adams--Novikov filtration, and contain occasional omissions. For this reason, we will reference \cite{brunerrognes} for specific calculations. 


\subsection{Defining the subgroup $\Free$ of $\pi_\ast\Tmf$}\label{achoiceofgenerators}

One often defines the elements in $\pi_\ast\Tmf$ by choosing a representative from the $E_2$-page of the descent spectral sequence of \cite{konter} or \cite{smfcomputation}. As with any spectral sequence though, we only know these elements are well-defined up to higher filtration. In this section, we define the subgroup $\Free\subseteq\pi_\ast\Tmf$, which in the reader's mind should be ``elements in $\pi_\ast\Tmf$ of lowest filtration in the descent spectral sequence'', but which we need to make precise below. A particular subtlety occurs at the prime $2$: we cannot explicitly define $\Free_d\subseteq\pi_d\Tmf_2$ for $d\equiv_{192}60, 156$ and $d>0$, a problem also encountered by Bruner--Rognes \cite[Rmk.9.24(4)]{brunerrognes}.\\

First, let us start with the following lemma.

\begin{lemma}\label{determiningelements}
Implicitly localise $\tmf$ at the prime $2$ and write $\mf_k=H^0(\M_{\Ell,\Z_{(2)}},\omega^{\otimes k})$ for the group of weight $k$ holomorphic modular forms over $\Z_{(2)}$. Then the following elements uniquely exist in $\pi_\ast\tmf$:
\begin{enumerate}
\item A class $c_4\in \pi_8\tmf$ which maps to the normalised Eisenstein series $c_4\in \mf_4$ and is $\bar{\kappa}$-torsion.
\item For each $k\in \{0,1,2,3,4,5,6\}$, a class $[c_4\Delta^{k+1}]\in \pi_{32+24k}\tmf$ which maps to $c_4\Delta^{k+1}\in \mf_{16+12k}$ and is $\bar{\kappa}$-torsion.
\end{enumerate}
\end{lemma}

\begin{proof}
In the first case, we can choose our $c_4$ to be $\widetilde{B}$ in the notation of \cite[Df.9.22 \& 9.50]{brunerrognes}. By \cite[Pr.9.40]{brunerrognes}, we see that $\bar{\kappa}c_4=0$. This class is determined by the fact that it detects $c_4$ and is $\bar{\kappa}$-torsion. Indeed, the ambiguity in our choice of $c_4$ lies in a factor of $\epsilon$, the image of $\epsilon\in \pi_8\Sph$. This ambiguity is solved by computing the Adams--Novikov spectral sequence (ANSS) for $\tmf$ of \cite[\textsection8]{bauer}, where we explicitly see that $\epsilon\bar{\kappa}\neq 0$ on the $E_2$-page.\\

In the second case, notice that the $k=1,5,6$ are uninteresting as there is no torsion class in these degrees, so the edge map in the Adams--Novikov spectral sequence is injective. Otherwise, we choose $[c_4\Delta^{k+1}]$ to be $\widetilde{B}_{k+1}$ in the notation of \cite[Df.9.22 \& 9.50]{brunerrognes}. By \cite[Lm.9.11]{brunerrognes}, we see that $\bar{\kappa}$ is $B=c_4+\epsilon$ power torsion, so \cite[Cor.9.55]{brunerrognes} states that $\widetilde{B}_{k+1}\bar{\kappa}=0$ for all $k$. Similar to the first case, these classes are uniquely defined by these properties. Indeed, the ambiguity of this choice is up to the higher filtration elements $\epsilon_k$ above $\widetilde{B}_{k+1}$ for $k\in \{0,3,4\}$ and $\bar{\kappa}$ above $\widetilde{B}_3$. An inspection of the ANSS $E_2$-page shows these torsion classes support nontrivial multiplication by $\bar{\kappa}$.
\end{proof}

We can now move onto our basis of $\Free$. For a modular form $f$ of weight $k$, we will write $f$ for an element in $\pi_{2k}\tmf$ which maps to $f$ under the edge map if such an element is uniquely determined by this fact or if this element is mentioned in \Cref{determiningelements}. 

\begin{notation}\label{elementsintmf}
The elements of $\Tors\subseteq\pi_\ast\tmf$ are simply the torsion elements, which can also be interpreted as elements in strictly positive filtration in the Adams--Novikov spectral sequence (ANSS)---this spectral sequence is called the \emph{elliptic spectral sequence} in \cite[\textsection7-8]{bauer}, which is identified with the desired ANSS in \cite[\textsection5]{tmfhomology} using the Gap Theorem of \cite{konter} or \cite{smfcomputation}. The elements of $\Free\subseteq \pi_\ast \tmf$ in nonnegative degree are then described in the following three cases:
\begin{itemize}
\item When $6$ is inverted, $\Free=\pi_\ast\tmf[\frac{1}{6}]$ as there is no torsion.
\item When localised at $3$, $\Free$ is multiplicatively generated by the classes:
\[c_4, c_6, [3\Delta], [c_4\Delta], [c_6\Delta], [3\Delta^2], [c_4\Delta^2], [c_6\Delta^2], \Delta^3\]
\item When localised at $2$, $\Free$ is multiplicatively generated by the classes
\[c_4, [2c_6], [8\Delta^{2i+1}], [4\Delta^{4j+2}], [2\Delta^4], [c_4\Delta^{k+1}], [2c_6\Delta^{k+1}], \Delta^8\]
for $i\in \{0,1,2,3\}$, $j\in\{0,1\}$, and $k\in\{0,1,2,3,4,5,6\}$, using \Cref{determiningelements} when necessary and where $[2c_6\Delta^{k+1}]$ is defined with additive indeterminacy $2\bar{\kappa}^3$ and $\eta\nu_6\epsilon$ for $k=1$ and $5$, respectively; see \cite[Rmk.9.24(4)]{brunerrognes}.
\end{itemize}
Define $\Free \subseteq \pi_\ast\Tmf$ in nonnegative degrees as the subset $\Free\subseteq \pi_\ast \tmf$ given above and in negative degrees as follows:
\begin{itemize}
\item When $6$ is inverted, then there is no torsion and the $\Z[\frac{1}{6}]$-module $\Free$ is generated by elements of the form $\{c_4^ic_6^j\Delta^k\}$, for $i\leq -1$, $j\in\{0,1\}$, and $k\leq -1$.; see \cite[Th.3.1]{konter}.
\item When localised at $3$, $\Free$ is generated by elements of the form
\[\{c_4^ic_6^j\Delta^{k-3l}\}, \{c_4^{-1}c_6\Delta^{-1-3l}\}, \{\frac{1}{3}c_4^{-1}c_6\Delta^{-2-3l}\}, \{\frac{1}{3}c_4^{-1}c_6\Delta^{-3(l+1)}\},\]
where $i\leq -1$, $j\in \{0,1\}$, $k\in \{-3,-2,-1\}$, $j+k<0$, and $l\geq 0$; see \cite[Th.4.1]{konter}. Unlike \cite{konter}, we have used the brackets $\{-\}$ to express torsion-free classes in negative degree to remind us that there is a degree shift that differs from the torsion-free classes in positive degree:
\[\{c_4^ic_6^j\Delta^k\}\in \pi_{8i+12j+24k-1}\Tmf_{(3)}.\]
The classes $\{c_4^{-m}\Delta^{-n}\}\in \pi_{-8n-24m-1}\Tmf_{(3)}$ for positive $n$ and $m$ such that $-8n-24m-1\equiv -49$ modulo $72$ above are not necessarily well-defined by their representative on the $E_2$-page of the descent spectral sequence, so we define them as the product of two well-defined elements $c_4\{c_4^{-m-1}\Delta^{-n}\}$.
\item When localised at $2$, $\Free$ is generated by elements of the form
\[\{c_4^i2c_6^j\Delta^{k-8l}\}, \{c_4^{-1}c_6\Delta^{k-8l}\}, \{c_4^{-1}2^{e_2(k+1)-2}c_6\Delta^{k-8l}\},\]
where $i<-1$, $j\in \{0,1\}$, $k\in \{-8,-7,\ldots, -1\}$, $l\geq 0$, and $e_2$ is the function which sends a nonzero integer $a$ the largest integer $b$ with $2^b|a$, and $e_2(0)=3$; see \cite[Th.5.3]{konter}. Similar to the $3$-local case, any potentially ambiguous elements can be defined as the product of either $c_4$ or $c_4^2$ with another well-defined element. For instance, we define the element $\{c_4^i2c_6^j\Delta^k\}$ in $\pi_q \Tmf_{(2)}$ for some negative $q=8i+12j+24k-1$, as the product of $c_4\{c_4^{i-1}2c_6^j\Delta^k\}$ if $q$ is congruent modulo $192$ to an element in the set
\[\{-37, -57, -61, -81, -97, -121, -133, -153, -157, -177\},\]
and as the product $c_4^2\{c_4^{i-2}2c_6^j\Delta^k\}$ if $q$ is congruent modulo $192$ to an element in the set
\[\{-49, -73, -145, -169\}.\]
\end{itemize}
\end{notation}


Let us now explicitly compare our generators to \cite{brunerrognes} (see Definitions 9.22, 9.50, and 13.13) in positive degrees: first at the prime $3$, then at the prime $2$---thank you to John Rognes for noticing our previous misreading of \emph{loc.\ cit}.

\begin{center}
\begin{tabular}{|c|c|c|c|c|c|c|c|c|c|c|} 
 \hline
$\Free$			&	$c_4$&	$c_6$&	$[3\Delta]$&	$[c_4\Delta]$&	$[c_6\Delta]$&	$[3\Delta^2]$&	$[c_4\Delta^2]$&	$[c_6\Delta^2]$&	$\Delta^3$	\\
 \hline
\cite{brunerrognes}	&	$B=B_0$&	$C=C_0$&	$D_1$&	$B_1$&		$C_1$&		$D_2$&		$B_2$&		$C_2$&		$H$	\\
 \hline
\end{tabular}
\end{center}

\begin{center}
\begin{tabular}{|c|c|c|c|c|c|c|c|c|c|} 
 \hline
$c_4$&	$[2c_6]$&	$[8\Delta^{2i+1}]$&	$[4\Delta^{4j+2}]$&	$[2\Delta^4]$&	$[c_4\Delta^{k+1}]$&	$[2c_6\Delta^{k+1}]$&	$[2c_6\Delta^{k'+1}]$	&	$\Delta^8$	\\
 \hline
$\widetilde{B}$&	$C$&	$D_{2i+1}$&	$D_{4j+2}$&		$D_4$&		$\widetilde{B}_{k+1}$&		$C_{k+1}$&		$C_{k'+1}+?$	&	$M$	\\
 \hline
\end{tabular}
\end{center}

Above we write $k\in \{0,2,3,4,6,7\}$ and $k'\in \{1,5\}$ and the question mark above indicates that those the objects $C_2$ and $C_6$ are only well-defined up to some additive indeterminacy.


\subsection{Anderson duality}\label{andersondualitysection}
To systematically study the negative homotopy groups of $\Tmf$, we will use the following form of duality.

\begin{mydef}
For an injective abelian group $J$, we write $I_J$ for the spectrum represented by the cohomology theory
\[\Sp\to \Ab_\ast\qquad X\mapsto \Hom_{\Ab_\ast}(\pi_{-\ast}X, J).\]
For a general abelian group $A$, we take an injective resolution of the form $0\to A\to J_1\to J_2$, which by functoriality, yields a morphism of spectra $I_{J_1}\to I_{J_2}$. The fibre of this morphism we denote by $I_A$, and for a spectrum $X$, we define the \emph{Anderson dual of $X$} to be the function spectrum $I_A X=F(X, I_A)$.
\end{mydef}

From the definition above one can calculate
\[\pi_\ast I_J X\simeq \Hom_{\Z}(\pi_{-\ast} X, J)\]
for an injective abelian group $J$. When $A$ is a general abelian group, we obtain the following functorial exact sequence of abelian groups for all $k\in \Z$
\begin{equation}\label{homotopyofad} 0\to \Ext_{\Z}^1(\pi_{-k-1}X, A)\to \pi_k I_A X\to \Hom_\Z(\pi_{-k}X, A)\to 0\end{equation}
which non-canonically splits when $A$ is a subring of $\Q$. More basic facts about Anderson duality, such as the fact that the natural map $X\to I_AI_A X$ is an equivalence when $X$ has finitely generated homotopy groups, can be found in \cite[\textsection6.6]{sagname}, under the guise of \emph{Grothendieck duality} in spectral algebraic geometry. Anderson duality is of interest to us as many of the spectra we will study in this article are \emph{Anderson self-dual}.

\begin{mydef}
Let $X$ be a spectrum and $A$ an abelian group. We say that $X$ is \emph{Anderson self-dual} if it comes equipped with an integer $d$ and an equivalence of spectra
\[\phi\colon X[d]\xrightarrow{\simeq} I_A X.\]
We also want to define a stricter form of self-duality for ring spectra. Let $R$ be an $\E_1$-ring with $\pi_0R\simeq A$ such that $\pi_{-d} R$ is a free $A$-module of rank one. We say an element $D\in \pi_{-d}R$ \emph{witnesses the Anderson self-duality of $R$} if the isomorphism $\phi_D\colon \pi_{-d} R\to A$ sending $D\mapsto 1$ which identifies $D$ as an $A$-module generator of $\pi_{-d}R$, lifts to an element $D^\vee\in \pi_d I_A R$ under the surjection of (\ref{homotopyofad}) whose representing map of left $R$-modules $D^\vee\colon R[d]\to I_A R$ is an equivalence.
\end{mydef}

\begin{example}\label{exampleofselfduality}
There are some famous examples of Anderson self-duality.
\begin{itemize}
\item The class $1\in \pi_0\KU$ witnesses the Anderson self-duality of $\KU$, ie,
\[1^\vee\colon \KU\xrightarrow{\simeq} I_\Z \KU\]
is an equivalence. This is originally due to Anderson \cite{andersonduality}, and is an immediate consequence of the fact that $\Hom_\Z(\pi_\ast\KU, \Z)$ is a free $\pi_\ast\KU$-module; see \cite[p.3]{drewvesna}.
\item The class $vu_\R^{-1}\in \pi_{-4}\KO$ witnesses the Anderson self-duality of $\KO$, ie,
\[(vu_\R^{-1})^\vee\colon \KO[4]\xrightarrow{\simeq} I_\Z \KO\]
is an equivalence. This result is also due to Anderson. An accessible modern proof with an eye towards spectral algebraic geometry can be found in \cite[Th.8.1]{drewvesna}.
\item The class $D=\{2c_4^{-1}c_6\Delta^{-1}\}\in \pi_{-21}\Tmf$ witnesses the Anderson self-duality of $\Tmf$, ie,
\[D^\vee\colon \Tmf[21]\xrightarrow{\simeq} I_\Z \Tmf\]
is an equivalence. The abstract duality result, meaning the existence of such an equivalence of $\Tmf$-modules above, is due to Stojanoska; see \cite[Th.13.1]{vesnaduality} for the case with 2 inverted and \cite{vesnaprimetwo} where it is announced in general; the $2$-primary case can also be found in \cite[Th.10.13]{brunerrognes}. Any such equivalence of $\Tmf$-modules is \emph{a posteriori} defined by a generator of $\pi_{21}I_\Z\Tmf\simeq \Z$, which we choose to be the above $D^\vee$, dual to $D$ using (\ref{homotopyofad}).
\end{itemize}
\end{example}

There are other examples for self-duality of topological modular forms \emph{with level structure}, as discussed for $\Tmf(2)$ in \cite[Th.9.1]{vesnaduality} and $\Tmf_1(m)$ in \cite[Th.5.14]{tmfwls}. Studying endomorphisms of Anderson self-dual spectra leads us to \emph{dual endomorphisms}.

\begin{mydef}\label{dualendomorphism}
Let $A$ be an abelian group, $X$ an Anderson self-dual spectrum, and $F\colon X\to X$ an endomorphism of $X$. Define the \emph{dual endomorphism} of $F$ as the composite
\[\widecheck{F}\colon X\xrightarrow{\phi, \simeq} (I_A X)[-d]\xrightarrow{(I_A F)[-d]} (I_A X)[-d]\xleftarrow{\phi, \simeq} X.\]
\end{mydef}

Given $A,X$, and $F$ from the above definition, then the functoriality of (\ref{homotopyofad}) yields the following commutative diagram of abelian groups with exact rows for all $k\in\Z$:
\begin{equation}\label{functorialityofhomotopyofad}\begin{tikzcd}
{0}\ar[r]		&	{\Ext_{\Z}^1(\pi_{-k-1-d}X, A)}\ar[r]\ar[d, "{\Ext^1_\Z(F, A)=F^\ast_1}"]	&	{\pi_k X}\ar[r]\ar[d, "{\widecheck{F}}"]	&	{\Hom_\Z(\pi_{-k-d}X, A)}\ar[d, "{\Hom_\Z(F, A)=F^\ast_0}"]\ar[r]		&	{0}	\\
{0}\ar[r]		&	{\Ext_{\Z}^1(\pi_{-k-1-d}X, A)}\ar[r]							&	{\pi_k X}\ar[r]					&	{\Hom_\Z(\pi_{-k-d}X, A)}\ar[r]								&	{0}
\end{tikzcd}\end{equation}
Our calculations of $\psi^k$ on $\Tmf$ in negative degrees will rest upon explicit calculations of $\widecheck{\psi}^k$ on positive homotopy groups and (\ref{functorialityofhomotopyofad}).\\

When working with $6$ inverted, there also exists a kind of algebro-geometric duality on $\M_\Ell$ called \emph{Serre duality}. The following can be found in \cite[\textsection A]{adddecompformf} using the well-known identification of $\M_{\Ell, \Z[\frac{1}{6}]}$ with the weighted projective line $\P_{\Z[\frac{1}{6}]}(4,6)$; see \cite[Ex.2.1]{adddecompformf}.

\begin{theorem}\label{rationalserreduality}
The dualising sheaf for $\M_{\Ell, \Z[\frac{1}{6}]}$ is $\omega^{-10}$. In particular, for any integer $k$ the natural cup product map
\[H^0(\M_{\Ell, \Z[\frac{1}{6}]}, \omega^k)\otimes H^1(\M_{\Ell, \Z[\frac{1}{6}]}, \omega^{-k-10})\to H^1(\M_{\Ell, \Z[\frac{1}{6}]}, \omega^{-10})\simeq \Z[\frac{1}{6}]\]
is a perfect pairing of $\Z[\frac{1}{6}]$-modules.
\end{theorem}

Let us note that the stack $\M_{\Ell}$ certainly has \textbf{no} Serre duality before inverting $6$, which can be seen through the cohomology calculations of $\omega^{\ast}$ over $\M_{\Ell}$ from \cite{konter}. 

\begin{remark}\label{rationaldualitygivesalittlehomotopy}
A simple consequence of the above theorem is that one can immediately see the $\E_\infty$-ring $\Tmf[\frac{1}{6}]$ is Anderson self-dual. Indeed, as discussed on \cite[p.8]{vesnaduality}, the Serre duality statement of \Cref{rationalserreduality}, the calculation of $H^\ast(\M_{\Ell,\Z[\frac{1}{6}]}, \omega^\ast)$ in \cite[\textsection3]{konter}, and a collapsing DSS, immediately implies the Anderson self-duality of $\Tmf[\frac{1}{6}]$ as in \Cref{exampleofselfduality}.
\end{remark}

When $6$ is inverted, dual endomorphisms on $\Tmf$ (defined using Anderson duality) can be computed directly using Serre duality.

\begin{lemma}\label{serredualityinaction}
Let $\P$ be a set of primes containing both $2$ and $3$ and implicitly localise everywhere away from $\P$. If $F\colon \Tmf\to \Tmf$ is a morphism of spectra, then one can compute $\widecheck{F}$ on $\pi_\ast\Tmf$ in negative degrees as the composite
\[\widecheck{F}\colon \pi_k\Tmf\simeq H^0(\M_\Ell, \omega^{-\frac{k+1}{2}-10})^\vee \xrightarrow{F^\vee} H^0(\M_\Ell, \omega^{-\frac{k+1}{2}-10})^\vee\simeq \pi_k\Tmf\]
and in nonnegative degrees as the composite
\[\widecheck{F}\colon \pi_k\Tmf\simeq H^1(\M_\Ell, \omega^{-\frac{k}{2}-10})^\vee \xrightarrow{F^\vee} H^1(\M_\Ell, \omega^{-\frac{k}{2}-10})^\vee\simeq \pi_k\Tmf,\]
where we have implicitly used the Serre duality isomorphism.
\end{lemma}

\begin{proof}
This follows immediately from the definitions, as in this case, the Anderson duality equivalence comes directly from Serre duality; see \Cref{rationaldualitygivesalittlehomotopy}.
\end{proof}


\subsection{Proof of \Cref{thb}}\label{generaloperations}

To prove \Cref{thb} we will use the following lemmata, the first helping us to calculate inside $\Free$ and the second to help us with $\Tors$.

\begin{lemma}\label{filtrationzeroandhure}
Let $R$ be an algebra in $\h\Sp$ and $A$ an $R$-algebra in $\h\Sp$. Suppose we have a decomposition of $\pi_\ast A$ given by $\Tors\oplus\Free$, where the elements of $\Tors$ are precisely the $\pi_0R$-torsion elements of $\pi_\ast A$. Fix some $\pi_k A$. Suppose that for each $y\in \Tors\subseteq \pi_k A$, there is a $z$ in the image of the unit $\pi_\ast R\to\pi_\ast A$ such that $zx=0$ for all $x\in \Free\subseteq \pi_\ast A$ and the map of $\pi_0 R$-modules
\begin{equation}\label{multiplicationinjective}\pi_kA\supseteq \langle y\rangle \xrightarrow{z\cdot} \langle zy\rangle \subseteq \pi_{k+|z|} A\end{equation}
is injective. Then for every $R$-module map $F\colon A\to A$, the induced map on homotopy groups $F\colon \pi_k A\to \pi_k A$ preserves the decomposition $\Tors\oplus\Free$.
\end{lemma}

\begin{proof}
Clearly $F(\Tors)\subseteq \Tors$ as $F$ is $R$-linear. Take an $x\in \Free$ and write $F(x)=x'+y$ where $x'\in \Free$ and $y\in\Tors$ using the decomposition above. The hypotheses then lead us to the equalities
\[0=F(zx)=zF(x)=z(x'+y)=zy\]
where the second equality follows from the $R$-linearity of $F$. The injectivity of (\ref{multiplicationinjective}) leads us to the conclusion that $y=0$, and we are done.
\end{proof}

\begin{lemma}\label{torsionistrivial}
Let $p$ be a prime, $x$ be a homogeneous element of $\Tors\subseteq\pi_\ast\tmf_p$, and $F\colon \tmf_p\to\tmf_p$ a morphism of spectra. Furthermore, if $p=2$, suppose that on rational homotopy groups we have the equality $F(c_4^m\Delta^l)=\lambda_{m,l} F(1) c_4^m\Delta^l$ for all $m\geq 1$ and $l\geq 0$, where $\lambda_{m,l}$ is an integer congruent to $1$ modulo $8$. Then we have the equality
\[F(x)=xF(1)\qquad \in \pi_\ast\tmf_p.\]
\end{lemma}

It will become clear during the proof that the above hypotheses can be somewhat weakened, but we will not need any generalisation in this article.\\

The following proof is quite long and relies on a case-by-case analysis of $\pi_\ast\tmf$.

\begin{proof}
Let us start by considering two purely formal cases:

\begin{enumerate}
\item Suppose $x$ is in the image of the Hurewicz morphism $\pi_\ast\Sph\to \pi_\ast\tmf$. These classes are displayed in colour in \cite[\textsection 13, p.2-4]{tmfbook} as conjectured by Mahowald and recently proven at the prime 2 in \cite{hurewicztmf} and at the prime 3 in \cite{tmfthree}; these facts can also be found in \cite[\textsection11.11 \& \textsection13.7]{brunerrognes}, respectively. In this case, as $F$ is $\Sph$-linear (all maps of spectra are) we obtain an equality:
\[F(x)=xF(1)\qquad \in \pi_{|x|}\tmf\]
\item Suppose that there exists an element $y\in\pi_\ast\tmf$ in the Hurewicz image such that $xy$ lies in the Hurewicz image and that the multiplication-by-$y$ map
\begin{equation}\label{injectivemultplication}\cdot y\colon \pi_{|x|}\tmf \to \pi_{|xy|}\tmf\end{equation}
is injective. In this case, we have the equalities
\[xyF(1)=F(xy)=F(x)y\]
which using the assumption that (\ref{injectivemultplication}) is injective, implies that $F(x)=xF(1)$. 
\end{enumerate}

These first two cases cover all of the torsion at the prime 3, so let us now focus on the prime 2. Consider the family of elements of the form
\[\eta^ic_4^j\Delta^k\qquad i\in\{1,2\} \quad j,k\geq 1\]
where we have temporarily foregone the use of brackets. For these elements, we immediately obtain the equality
\[F(\eta^ic_4^j\Delta^k)=\eta^i F(c_4^j\Delta^k)\in \pi_{i+8j+24k}\tmf_{2}\]
using $\Sph$-linearity. Moreover, we claim to have the equality
\begin{equation}\label{goalequiationone}F(c_4^j\Delta^k)=\lambda_{j,k} c_4^j\Delta^k F(1)\in \pi_{8j+24k}\tmf_{2}.\end{equation}
Indeed, this is the na\"{i}ve calculation from the zero line in the $E_2$-page of the ANSS and the fact that this zero line injects into the rational homotopy groups. Hence we need to check that $F$ preserves $\Free_d$ in these degrees $d=8j+24k$, ie, that $F$ does not jump filtration in these degrees. The only nontrivial cases to check, so those $d$ where $\Tors_d\neq 0$, are those $d$ congruent modulo $192$ to a number in the set
\begin{equation}\label{problematicdegrees}	\{8,32,40,80,104,128,136\}.	\end{equation}
In the above cases, we have potential torsion classes
\[\epsilon,\quad q=\epsilon_1,\quad\bar{\kappa}^2,\quad \bar{\kappa}^4,\quad [\epsilon\Delta^4]=\epsilon_4, \quad[q\Delta^4]=\epsilon_5,\quad \bar{\kappa}^2[2\Delta^4]=\bar{\kappa}^2 D_4\]
where we have used the notation of \cite{bauer} on the left and \cite[Df.9.22]{brunerrognes} on the right. We want to apply \Cref{filtrationzeroandhure} with $z=\bar{\kappa}$. The fact that $yz\neq0$ for all $y$ in the set (\ref{problematicdegrees}) above follows from \cite[Pr.9.41]{brunerrognes} and the fact that $\bar{\kappa}^5\neq 0$. This allows us to use \Cref{filtrationzeroandhure} to conclude the equality (\ref{goalequiationone}) for those (\ref{problematicdegrees}). The fact that for all remaining $d$ and all $x\in \Free_d$, we have $x\bar{\kappa}=0$ follows from the facts that $[c_4\Delta^{k}]\bar{\kappa}=0$ for all $0\leq k\leq 7$ as shown in the proof of \Cref{determiningelements}, and the fact that all other classes in these degrees are $c_4$ and $\Delta^8$ multiplies of these classes. To summarise this argument thus far, we have the equalities
\[F(\eta^ic_4^j\Delta^k)=\eta^i F(c_4^j\Delta^k)=\eta^i \lambda_{j,k}c_4^j\Delta^k F(1)=\eta^i c_4^j \Delta^k F(1)\]
the latter coming from our hypothesis. We claim it suffices to now consider the two families of elements
\begin{equation}\label{forthcase} 
[\eta\Delta]=\eta_1,\quad [\eta\Delta]^2=\eta_1^2,\quad [\eta\Delta]^3=\eta_1^3,\quad [\eta\Delta^4]=\eta_4,\quad  [\eta^2\Delta^5]=\eta_1\eta_4
\end{equation}
\begin{equation}\label{twoextensioncases}
[2\nu\Delta]=\nu_1,\quad [\nu\Delta^2]=\nu_2,\quad [\nu\Delta^4]=\nu_4,\quad [2\nu\Delta^5]=\nu_5, \quad [\nu\Delta^6]=\nu_6
\end{equation}
up to $\Delta^8$-periodicity and multiplication by an element in the Hurewicz image---again we have used the notation from \cite{bauer} on the left and that from \cite{brunerrognes} on the right of the equalities. We will now detail an argument for the element $[\eta\Delta]$, and all other elements of the first family (\ref{forthcase}) follow similarly. Our map $F\colon \tmf_{2}\to \tmf_{2}$ of spectra induces a map of ANSSs. As our original map is $\Sph$-linear and the ANSS functor is lax-monoidal,\footnote{Indeed, the ANSS can be viewed as the cobar spectral sequence of the cosimplicial spectrum $X\otimes \MU^{\otimes (\bullet+1)}$, which induces a lax-monoidal functor by \cite[Lm.2.39]{achimthesis}. One could also use the more classical and direct argument found in \cite[Th.2.3.3]{greenbook}.} this induced map of spectral sequences is linear over the ANSS for the sphere $\Sph$. The class $[\eta\Delta]$ has $E_2$-representative $h_1\Delta$ as we can see on \cite[p.32]{bauer}, where $h_1$ is the image of the class of the same name in the ANSS for $\Sph$ induced by the unit map $\Sph\to\tmf$. The value of the map induced by $F$ on the $E_2$-page can then be calculated as
\[F_2(h_1\Delta)=h_1 F_2(\Delta)\in E_2^{1,25}\simeq \Z/2\Z\{h_1\Delta\}\oplus\Z/2\Z\{h_1c_4^3\}.\]
Everything on the zero line of this $E_2$-page is torsion-free, so this line maps injectively into its rationalisation. Rationally, however, the ANSS for $\tmf$ collapses on the $E_2$-page, so $F_2(\Delta)$ may be calculated as $F(\Delta)$ inside $\pi_{24}\tmf_\Q$. From our hypotheses, we see that $F_2(\Delta)=\lambda\Delta F(1)$, where $\lambda$ is odd. This immediately yields the equality $F_2(h_1\Delta)=h_1\Delta F(1)$ inside $E_2^{1,25}$. This equality also exists on the $E_\infty$-page $E_\infty^{1,25}$, and due to the lack of classes of higher filtration in the 25-stem, we obtain this equality in $\pi_{25}\tmf_{2}$, meaning $F([\eta\Delta)=[\eta\Delta]F(1)$. The argument works similarly for the other elements in the first family (\ref{forthcase}) as there are no classes of higher Adams--Novikov filtration in each degree considered above, the $E_2$-representative $h_1$ for $\eta$ comes from the ANSS for $\Sph$ and is $2$-torsion.\\

The case for the elements in the family (\ref{twoextensioncases}) follows similarly, except we need to be careful about the exotic $2$-extensions supported by these classes. In other words, it is no longer clear that our argument on the $E_2$-page carries over. To fix this, we will work with the synthetic spectrum $\nu\tmf/\tau^4$, which acts as an intermediary between the $E_2$-page and $E_\infty$-page. This remedy was suggested to us by an anonymous referee, who we heartily thank---another thank you to Christian Carrick for helping us out with some details below.\\

Let us consider the argument for $[\nu\Delta^2]=\nu_2$---the other cases follow with the obvious changes. Consider the $\BP$-synthetic category $\Syn$ of \cite{syntheticspectra} at the prime $2$,\footnote{This is \textbf{not} the \emph{even} $\BP$-synthetic category. Inside $\Syn$, an element in $\pi_{\ast,\ast}\nu X$ is $\tau^{r-1}$-torsion if and only if it's hit by a $d_r$-differential in the $\BP$-baseed ANSS for $X$.} and in particular the $C(\tau^4)$-module internal to this category $X=\nu(\tmf)/\tau^4$, where $\nu\colon \Sp\to \Syn$ is the synthetic analogue functor. The rest of this proof also goes through without major changes with the synthetic spectrum $\mathrm{Smf}$ of \cite[Th.C]{osyn} and its $\sigma$-spectral sequence computed in \cite{smfcomputation}. We calculate $\pi_{51,\ast}\nu(\tmf)$ to be the $\Z_{(2)}[\tau]$-module
\[\Z/8\Z[\tau]\{[\nu\Delta^2]\}\oplus (\Z/2\Z[\tau]/\tau^2)\{h_0^3c_4^6,h_0^3c_4^3\Delta\}\oplus V_{51}\]
where $V_{51}$ is all $2$- and $\tau^2$-torsion and comes from elements in filtration 7 or higher in the ANSS. Similarly, we can calculate $\pi_{d,\ast}X$ for $d=48$ and $51$ as the following $\Z_{(2)}[\tau]/\tau^4$-modules:
\[\pi_{48,\ast}X\simeq (\Z_{(2)}[\tau]/\tau^4)\{\Delta^2, c_4^3\Delta,c_4^6\}\oplus W_{48}\]
\[\pi_{51,\ast}X\simeq (\Z/8\Z[\tau]/\tau^4)\{\nu\Delta^2\}\oplus (\Z/2\Z[\tau]/\tau^2)\{h_0^3c_4^6,h_0^3c_4^3\Delta\}\oplus W_{51}\]
The $\Z_{(2)}[\tau]/\tau^4$-modules $W_{48}$ and $W_{51}$ above are both $2$- and $\tau^2$-torsion and come from elements in filtration $4$ or higher. Notice that we still have an element $\nu\Delta^2$ which is strictly $8$-torsion in $\pi_{51,\ast}X$ as the exotic multiplication by $2$ jumps only $2$-filtrations in the $E_\infty$-page of the classical ANSS for $\tmf$, hence it is detected in $\nu(\tmf)/\tau^4$.\\

Our assumption about the effect of $F$ on rational homotopy groups implies that $F(\Delta^2)\equiv \Delta^2 F(1)$ modulo $8$ on the $E_2$-page of the classical ANSS for $\tmf$. The naturality of the $\tau$-Bockstein spectral sequence for $X=\nu(\tmf)/\tau^4$ implies that the map induced by $F$ on $\pi_{48,\ast}X$ sends $\Delta$ to itself modulo $8$ and elements in $W_{48}$. As there exists a lift of $\nu\in \pi_3\Sph$ inside $\pi_{3,\ast}\nu\Sph$, then by $\nu\Sph$-linearity we see that the map induced by $F$ sends the element $\nu\Delta$ to itself modulo elements in $\nu W_{48}\subseteq W_{51}$.\\

The canonical quotient map $\nu(\tmf)\to X$ sends the $[\nu\Delta^2]$ generating a $\Z/8\Z[\tau]$ to the $8$-torsion class $\nu\Delta^2$ inside $\pi_{51,\ast}X$. The naturality of $-\otimes C(\tau^4)$ and our calculation above shows that the map induced by $F$ on $\pi_{51,\ast}\nu(\tmf)$ sends $[\nu\Delta^2]$ to itself modulo $V_{51}$ and $\tau$-torsion. Inverting $\tau$, using the fact that $V_{51}$ is all $\tau^4$-torsion, we see that $F([\nu\Delta^2])=[\nu\Delta^2]$ inside $\pi_{51}\tmf$.\\

For other elements in (\ref{twoextensioncases}) the above argument runs through with the evident changes without any surprises.
\end{proof}

We are now in a position to prove \Cref{thb}.

\begin{proof}[Proof of \Cref{thb}]
First, we start with elements $x\in \Free\subseteq \pi_\ast\Tmf_p$ in nonnegative degrees. As the operations $\psi^k$ are multiplicative, it suffices to calculate $\psi^k$ on these generators of $\Free$. If $x$ lies in a degree with no torsion elements, then our calculation on the $E_2$-page of the descent spectral sequence (DSS) holds, and we are done. If there is torsion, in this degree, we have to make another argument. Checking our definition of $\Free$ in nonnegative degree and the homotopy groups of $\Tmf_p$, we first notice that at the prime $3$, there are no generators of $\Free$ in nonnegative degree with a nonzero torsion class also in that degree, so we focus on the case of $p=2$. At this prime, the only problematic nonnegative degrees lie in the following list congruent $192$:
\[8, 32, 60, 80, 104, 128, 156\]
In the cases other than $60$ and $156$, our $E_2$-page calculation yields the calculation on homotopy groups using \Cref{filtrationzeroandhure}, where $z=\bar{\kappa}$; this technique was already used in the proofs of \Cref{determiningelements} and \Cref{torsionistrivial}. In degrees $d$ congruent to $60$ and $156$ the group $\Free_d$ is not well-defined, hence the exception in these degrees.\footnote{The failure to find well-defined elements in these degrees is closely related to the failure of \Cref{filtrationzeroandhure} in these degrees as well, as the torsion classes in these degrees do not support any interesting multiplication.} To summarise, for all $x\in \Free_d$ for nonnegative $d$, the $E_2$-page calculations holds and we obtain $\psi^k(x)=k^{\frac{d}{2}}x$.\\

Suppose now that $x\in \Free$ has negative degree. Looking at our generators of $\Free$ as an abelian group from \Cref{elementsintmf} in negative degrees, we notice that these generators are either in a degree with no torsion or defined as the product of such a class with $c_4$ or $c_4^2$. From this observation, it suffices to calculate $\psi^k$ on $x\in \Free$ in negative degrees where there is no torsion, hence we may invert $p$ and work inside $\pi_\ast\Tmf_p[\frac{1}{p}]$. In this case, we have to compute the morphism
\[\psi^k\colon H^1(\M_{\Ell, \Q_p}, \omega^d)\to H^1(\M_{\Ell, \Q_p}, \omega^d)\]
for all $d<0$. This we can do with a calculation of the cohomology of the stack with graded structure sheaf $(\M_{\Ell, \Q_p}, \omega^\ast)$, which is equivalent to the weighted projective line $\P_{\Q_p}(4,6)$; see \cite[Ex.2.1]{adddecompformf}. In this case, we can use the fact that the groups $H^\ast(\P_{\Q_p}(4,6), \omega^\ast)$ are isomorphic to the groups $H^\ast(\widetilde{P}(4,6), \O)$, where $(\widetilde{P}(4,6),\O)$ is $(\Spec A-\{0\}, \O)$, where $A=\Q_p[c_4,c_6]$, together with the $\G_m$-action given by the gradings $|c_4|=4$ and $|c_6|=6$. As discussed for $\M(2)$ in \cite[\textsection7]{vesnaduality}, one can use the long exact sequence on cohomology induced by the expression $\widetilde{P}(4,6)\subseteq \Spec A \supseteq \{0\}$ \cite[Exercise III.2.3]{hart}, and the fact that $R\Gamma_{\{0\}}(\Spec A, \O)$ can be computed via the Koszul complex
\[A\to A[\frac{1}{c_4}]\times A[\frac{1}{c_6}]\to A[\frac{1}{c_4c_6}]\]
we obtain the following exact sequence
\[0\to A\to H^0(\widetilde{P}(4,6), \O)\to 0\to0\to H^1(\widetilde{P}(4,6), \O)\to A/(c_4^\infty, c_6^\infty)\to 0\]
Using this, we can explicitly calculate $\psi^k$ on $H^1(\widetilde{P}(4,6), \O)\simeq A/(c_4^\infty, c_6^\infty)$ as
\[\psi^k(\frac{1}{c_4^ic_6^j})=k^{-4i-6j}\frac{1}{c_4^ic_6^j}\]
where $\frac{1}{c_4^ic_6^j}$ represents a class in $\pi_\ast\Tmf_p[\frac{1}{p}]$ of topological degree $-8i-12j-1$. This yields the desired result.\\

Let us first consider a torsion element $x\in \Tors\subseteq \pi_\ast\Tmf_p$ and implicitly complete at $p$ for the rest of this proof. It suffices to consider the prime $p=2$ or $p=3$, otherwise $\Tors=0$. If $x$ has nonnegative degree, then we can immediately apply \Cref{torsionistrivial}, and we are done. Indeed, the hypotheses of that proposition apply as we already know $\psi^k(c_4^m\Delta^l)=k^{8m+12l}c_4^m\Delta^l$ and $k^{8m+12l}$ is congruent to $1$ modulo $8$ using Euler's theorem, for $k\in \Z_2^\times$. If $x$ is an element of $\Tors$ of negative degree, then we will consider (\ref{functorialityofhomotopyofad}) for $\Tmf$, which yields the commutative diagram of abelian groups for every integer $d$
\begin{equation}\label{tmfanddual}\begin{tikzcd}
{0}\ar[r]		&	{\Ext_{\Z}^1(\pi_{-d-22}\Tmf)}\ar[r]\ar[d, "{(\psi^k)_1^\ast}"]	&	{\pi_d \Tmf}\ar[r]\ar[d, "{\widecheck{\psi}^k}"]	&	{\Hom_\Z(\pi_{-d-21}\Tmf)}\ar[d, "{(\psi^k)_0^\ast}"]\ar[r]		&	{0}	\\
{0}\ar[r]		&	{\Ext_{\Z}^1(\pi_{-d-22}\Tmf)}\ar[r]	&	{\pi_d \Tmf}\ar[r]	&	{\Hom_\Z(\pi_{-d-21}\Tmf)}\ar[r]		&	{0}
\end{tikzcd}\end{equation}
where all Ext- and Hom-groups above have $\Z$ as a codomain and the dual operation $\widecheck{\psi}^k$ is defined in \Cref{dualendomorphism}. As $\psi^k$ induces a map of abelian groups on homotopy groups, we can then detect the effect of $\psi^k$ on $\Tors\subseteq\pi_\ast\Tmf$ by the effect of $(\psi^k)^\ast_1$ on the above Ext-groups. We want to use the Anderson self-duality diagram (\ref{tmfanddual}) to turn the computations of $\psi^k$ into computations of the dual operation $\widecheck{\psi}^k$ of \Cref{dualendomorphism}. In particular, we are reduced to compute the effect of $\widecheck{\psi}^k$ on elements in $\Tors$ of nonnegative degree, for which we would like to use \Cref{torsionistrivial}, again. This first requires us to calculate $\widecheck{\psi}^k(c_4^m\Delta^l)$ for $k\in \Z_2^\times$ after inverting $2$. Using \Cref{serredualityinaction} and the above calculations of $\psi^k$ to obtain the rational calculation
\[\widecheck{\psi}^k(c_4^m\Delta^l)=k^{-10-8m-12l}c_4^m\Delta^l.\]
We now use the fact that $(\Z/8\Z)^\times\simeq (\Z/2\Z)^2$ to see that for $k\in \Z_2^\times$, $k^{-10-8m-12l}$ is congruent to 1 modulo $8$. From this we see that \Cref{torsionistrivial} applies, which shows that for torsion elements $x$ in nonnegative degree, $\widecheck{\psi}^k(x)=x\widecheck{\psi}^k(1)=x k^{-10}=x$ as $k^{-10}\equiv_8 1$ for all $k\in \Z_2^\times$ and $k^{-10}\equiv_3 1$ for $k\in \Z_3^\times$, both of which are easily checked by hand. Using (\ref{tmfanddual}), we see that $\psi^k(x)=x$ for $x\in\Tors$ of negative degree, and we are done.
\end{proof}

The above proof shows that we can calculate $\widecheck{\psi}^k$ on $\pi_\ast\Tmf$ in certain degrees.

\begin{prop}\label{dualoperations}
Let $p$ be a prime and $k\in \Z_p^\times$ be a $p$-adic unit. Then the effect of the dual Adams operation $\widecheck{\psi}^k$ on $\pi_\ast \Tmf_p$ is given by
\[\widecheck{\psi}^k(x)=\begin{cases}	

x 					&	x\in \Tors	\\
k^{-10-\ceil{\frac{|x|}{2}}}x 	&	x\in \Free
\end{cases}\]
\textbf{unless} for $x\in\Free$ lies in the following degrees at the following primes:
\begin{itemize}
\item $p=3$ and $|x|=72r+40$ for some $r\geq 0$, then the answer holds modulo $\be^4\Delta^{3r}$.
\item $p=3$ and $|x|=72r-49$ for some $r<0$, then the answer holds modulo $\langle\al\be^2\Delta^{3(r-1)}\rangle$.
\item $p=2$ and $|x|=192 r+d$ for some $r\geq 0$ and $d$ in the set $\{20,60, 68, 100, 116, 156, 164\}$, then the answer holds modulo
\[4\bar{\kappa}\Delta^{8r}, 2\bar{\kappa}^3\Delta^{8r}, \kappa\nu[\nu\Delta^2]\Delta^{8r}, \bar{\kappa}^5\Delta^{8r}, 2\bar{\kappa}[2\Delta^4]\Delta^{8r}, \nu^3[\nu\Delta^6]\Delta^{8r}, \nu\kappa[\nu\Delta^6]\Delta^{8r},\]
respectively.
\item $p=2$, $|x|\geq 0$, and $|x|\equiv d$ modulo $192$, where $d$ is an element of
\[\{-49, -61, -73, -97, -121, -145, -157, -169\}\]
where the result holds modulo torsion.
\end{itemize}
\end{prop}

Our proof will follow the outline of the proof of \Cref{thb}, the only difference being that the operations $\psi^k$ are multiplicative and the $\widecheck{\psi}^k$ are not. A similar style of proof can be used to compute the effect of other operators on $\Tmf$ and $\TMF$ such as the Hecke operators of \cite{heckeontmf}.

\begin{proof}
If $x\in \Free_d$ lives in $\pi_d\tmf$ with no torsion, then the desired result follows by inverting $p$ and applying \Cref{serredualityinaction} and \Cref{thb}. If $\pi_d\tmf$ contains some torsion, then we want to apply \Cref{filtrationzeroandhure}. At the prime $3$, the only nonnegative degrees where we have problems are $d\equiv_{72}20,40$, the first is dealt with using \Cref{filtrationzeroandhure} with $z=\be$ and the latter case is an exception. Similarly for negative degrees; see \cite[\textsection4]{konter}. At the prime $2$, the problematic nonnegative degrees lie in the set
\[\{8,20,28,32,40,52,60,68,80,100,104,116,124,128,136,148,156,164\}\]
modulo $192$. All of the cases above can be dealt with using \Cref{filtrationzeroandhure} with $z=\kappa$, except for $d=20,40,60,80$ where we use $z=\bar{\kappa}$, and the exceptional cases. For the negative degrees, a similar we are reduced to degrees in the set
\[\{-37, -49, -57, -61, -73, -81, -97, -121, -133, -145, -153, -157, -169, -177\}\]
modulo $192$, which are again dealt with using \Cref{filtrationzeroandhure} or left as an exception. For $x$ in $\Tors$ in nonnegative degree, we can apply \Cref{torsionistrivial}. For $x$ in $\Tors$ in negative degree, we can look at (\ref{tmfanddual}) and use our calculations from \Cref{thb}.
\end{proof}

This extra calculation of $\widecheck{\psi}^k$ above suggests the following conjecture regarding the relation between endomorphisms and dual endomorphisms.

\begin{conjecture}\label{andersondualconjecture}
Let $R$ be an $\E_1$-ring spectrum and write $A=\pi_0 R$. Suppose that there is a class $D\in \pi_{-d} R$ such that $D$ witnesses the Anderson self-duality of $R$. Then, for any endomorphism $F\colon R\to R$ of algebra objects in $\h\Sp$ such that $F(D)=\lambda D$ for some $\lambda\in A$, the composites $F\circ \widecheck{F}$ and $\widecheck{F}\circ F$ are equivalent to multiplication by $\lambda$ on $\pi_\ast R$.
\end{conjecture}

An optimist might speculate that these potential equalities can perhaps be lifted to homotopies of morphisms of spectra.\\

This conjecture holds in the following cases:

\begin{itemize}
\item For $\KU_p$ and $\psi^k$ for $k\in \Z_p^\times$, one has $D=1$ and $\lambda=1$. In this case, the above conjecture can be checked using (\ref{functorialityofhomotopyofad}).
\item For $\KO_p$ and $\psi^k$ for $k\in \Z_p^\times$, one has $D=vu_\R^{-1}$ and $\lambda=k^{-2}$. In this case, the above conjecture can be checked using (\ref{functorialityofhomotopyofad}) again. Furthermore, Heard--Stojanoska verified that at the prime $2$ there is a homotopy between $\widecheck{\psi}^l$ and the $(-2)$-fold suspension of $\psi^{1/l}$, where $l$ is a topological generator of $\Z_2^\times/\{\pm 1\}$; see \cite[Lm.9.2]{drewvesna}.
\item For $\Tmf_p$ and $\psi^k$ for $k\in \Z_p^\times$, one has $D=\{2c_4^{-1}c_6\Delta^{-1}\}$ and $\lambda=k^{-10}$. In this case, the above conjecture can be checked (in some degrees) using \Cref{thb} and \Cref{dualoperations}.
\end{itemize}

\begin{remark}
Let us note a possible counter-example if we do not assume $F$ is multiplicative, as mentioned to us by Lennart Meier. Consider $F=\id+\psi^{-1}$ as an endomorphism of $\KU$. Then $\lambda=2$, however $F(u)=u-u=0$ on the usual generator $u\in \pi_2 \KU$, so \Cref{andersondualconjecture} cannot possibly hold in this case.
\end{remark}


\section{Applications}\label{applicationssection}

Our goal of this applications section is to show how one can easily manipulate the Adams operations on $\Tmf$ from \Cref{tha,thaversiontwo} as one does Adams operations on topological $K$-theory. In \Cref{adamssummandapplication}, we construct a connective height 2 Adams summand $\u$. That is, for each prime $p$ we define $\u=\tmf^{h\F_p^\times}$ using \Cref{thaversiontwo} such that $\u$ only has homotopy groups in nonnegative degrees divisible by the order of $v_1$, so divisible by $2(p-1)$. For $p=5$ the homotopy groups of this $\E_\infty$-ring $\u$ appear (meaning are isomorphic as a graded ring) to be of the form
\[\pi_\ast \u\simeq \Z_5[v_1, \sqrt{v_2}]\simeq \pi_\ast \BP\langle 2\rangle [\sqrt{v_2}]\]
which suggests this $\u$ is quite close to an $\E_\infty$-form of $\BP\langle 2\rangle$---similar observations also hold at the primes $7$ and $11$. We then prove \Cref{thd}, which states that $\tmf_p$ splits as a sum of shifts of $\u$ if and only if $p-1$ divides $12$, but when we invert $\Delta$ we always obtain the desired splitting. In \Cref{previouslydiscussedsplittingconjecturesection}, we conjecture that for primes $p$ such that $p-1$ does \textbf{not} divide $12$, there is a cofibre sequence involving a sum of shifts of $\u$, $\tmf_p$, and a sum of shifts of height $1$ Adams summands $\ell$. In \Cref{imageofjapplication}, we construct height 2 image-of-$J$ spectra $\s$ with maps $\s\to \j$ to the classical height 1 image-of-$J$ spectra. In particular, the fact that this map $\s\to \j$ is surjective on homotopy groups and the classical fact that $\Sph\to \j$ is split surjective on homotopy groups implies that $\Sph\to \s$ detects all of the image-of-$J$ elements inside $\pi_\ast\Sph$; see \Cref{the}. We hope that further refinements of $\s$ will bring us closer to a spectrum capturing height 1 and height 2 information, such as Behrens $Q(N)$ spectra do at large primes \cite{dividedbetafamily}, and that these spectra $\s$ might lend themselves to computations with an $\F_p$-based Adams spectral sequence.


\subsection{Connective height $2$ Adams summands and \Cref{thd}}\label{adamssummandapplication}

By \Cref{periodicitmfmotherufkcer}, we see that $\KU_p$ and $\TMF_p$ both have $p$-adic Adams operations $\psi^k$ for each $k\in\Z_p^\times$. When $p$ is odd, then $\Z_p^\times$ has a maximal finite subgroup $\F_p^\times$. This implies that both $\KU_p$ and $\TMF_p$ have $\E_\infty$-actions of the group $\F_p^\times$, which by a theorem of Gau{\ss} is isomorphic to the cyclic group of order $p-1$. A classical construction in homotopy theory is the \emph{Adams summand} $\KU_p^{h\F_p^\times}$, usually denoted by $L$, with connective cover $\ell$. Both $L$ and $\ell$ have simple homotopy groups as we are working with $p$-complete spectra and the group $\F_p^\times$ has order prime to $p$. In particular, we have isomorphisms
\[\pi_\ast\ell\simeq \Z_p[v_1]\qquad \pi_\ast L\simeq \Z_p[v_1^{\pm}]\]
where $v_1=u^{p-1}$ is the first Hasse invariant from chromatic homotopy theory and $\pi_\ast\KU\simeq \Z[u^{\pm}]$. When written like this, it is clear that $\ell$ is an $\E_\infty$-form of $p$-complete $\BP\langle 1\rangle$. These $\E_\infty$-rings $L$ and $\ell$ are summands of $\KU_p$ and $\ku_p$, respectively, associated with the idempotent map
\begin{equation}\label{idempotentusedtosplit}\frac{1}{p-1}\sum_{k\in \F_p^\times}\psi^k\end{equation}
revealing why they are called Adams \emph{summands}. In fact, more is true, as one can easily check that the canonical maps of $\E_\infty$-rings $L\to \KU_p$ and $\ell\to \ku_p$ recognise the codomain as a \emph{quasi-free}\footnote{Recall from \cite[Df.7.2.1.16]{haname} that for an $\E_\infty$-ring $R$ and an $R$-module $M$, we say $M$ is \emph{quasi-free} if there exists an equivalence $M\simeq\bigoplus_\al R[n_\al]$, and $M$ is \emph{free} if all of the $n_\al$ can be taken to be zero.} module over the source of rank $p-1$. Given we have the same $p$-adic Adams operations on $\Tmf_p$, we would like to explore the above ideas at the height two---the results are not what one might immediately expect; see \Cref{heighttwoadamssummandtheorm}. For an odd prime $p$, recall the $\F_p^\times$ action on the $\E_\infty$-rings $\TMF_p$ and $\tmf_p$ given by \Cref{periodicitmfmotherufkcer} and \Cref{thaversiontwo}, respectively. 

\begin{mydef}\label{heighttwoadamssummands}
For an odd prime $p$, define the $\E_\infty$-rings $\u=\tmf_p^{h\F_p^\times}$ and $\U=\TMF_p^{h\F_p^\times}$ and call them \emph{height two Adams summands}. For $p=2$ we set $\u=\tmf_2$ and $\U=\TMF_2$. By \Cref{thaversiontwo}, the natural map $\tmf_p\to \TMF_p$ is $\F_p^\times$-equivariant, factors through a map of $\E_\infty$-rings $\u\to \U$, and $\tmf_p\to \ku_p$ factors through a map of $\E_\infty$-rings $\u\to \ell$. In other words, we have the following commutative diagram of $\E_\infty$-rings
\begin{equation}
\begin{tikzcd}
	\u &&& \tmf_p \\
	& \U & \TMF_p \\
	& L & \KU_p \\
	\ell &&& \ku_p
	\arrow[from=1-1, to=1-4]
	\arrow[from=1-1, to=2-2]
	\arrow[from=1-1, to=4-1]
	\arrow[from=1-4, to=2-3]
	\arrow[from=1-4, to=4-4]
	\arrow[from=2-2, to=2-3]
	\arrow[from=3-2, to=3-3]
	\arrow[from=4-1, to=3-2]
	\arrow[from=4-1, to=4-4]
	\arrow[from=4-4, to=3-3]
\end{tikzcd}
\end{equation}
Recall that the map $f\colon \Tmf\to \KU$ does not extend over $\TMF$ as $f(\Delta^{24})=0$, so we also do not expect to see vertical maps between $\U$ and $L$ above.
\end{mydef}

We choose the names $\u$ and $\U$ as $\u$ is to $\tmf_p$ as $\ell$ is to $\ku_p$---we are open to other conventions. The homotopy groups of $\u$ and $\U$ as still simple to write down if $p\geq 5$:
\[\pi_\ast \u\simeq \left(\Z_p[x,y]\right)^{\F_p^\times}\simeq \Z_p\left\{x^i y^j | i,j\geq 0 \text{ such that }4i+6j \equiv_{p-1} 0\right\}\]
\[\pi_\ast \U\simeq \Z_p\left\{x^i y^j\Delta^{k} | i,j\geq 0, k\in\Z \text{ such that }4i+6j+12k \equiv_{p-1} 0\right\}\]
where $x=c_4$ has degree $8$, $y=c_6$ has degree $12$, and $\Delta=\frac{x^3-y^2}{1728}$. Both $\u$ and $\U$ are summands of $\tmf_p$ and $\TMF_p$, respectively, using the same idempotent (\ref{idempotentusedtosplit}) as the height one case. However, it is not true that the inclusion $\u\to \tmf_p$ witnesses the target as a quasi-free module over the source for all $p$ unlike the height one case.

\begin{theorem}\label{heighttwoadamssummandtheorm}
For every odd prime $p$ the map $\U\to \TMF_p$ recognises the codomain as a rank $\frac{p-1}{2}$ quasi-free module over the domain. The map $\u\to \tmf_p$ recognises the target as a rank $\frac{p-1}{2}$ quasi-free module if $p-1$ divides $12$ and for all other primes $\tmf_p$ is never a quasi-free $\u$-module.
\end{theorem}

Recall that $\u=\tmf_2$ and $\U=\TMF_2$ at the prime $2$, so we ignore this case above.\\

The proof of this theorem is rather elementary and consists of formal stable homotopy theory and some dimension formul{\ae} for spaces of (meromorphic) modular forms. We will write $\mf_\ast$ for the $p$-completion of $H^0(\M_{\Ell}, \omega^{\ast\geq 0})$ and $\MF_\ast$ for the $p$-completion of $H^0(\M_\Ell^\sm, \omega^\ast)$. Both of these cohomology rings are easy to calculate as in this case the $q$-expansion homomorphism into $\Z\llbracket q\rrbracket_p$ is injective; see \cite[Th.12.3.7]{diamondim}.

\begin{proof}
Let us start with the connective case---it is a little simpler. For $p=3$, the map $\u\to \tmf_3$ is an equivalence, as $\F_3^\times$ acts trivially on $\pi_\ast\tmf_3$ and the order of this group is invertible in $\pi_0\tmf_3\simeq \Z_3$ so the associated homotopy fixed point spectral sequence collapses. At $p=5$, we claim the map of $u$-modules
\[\u\oplus \u[12]\xrightarrow{1\oplus y}\tmf_5\]
defined by the elements $1,y\in \pi_\ast\tmf_5$, is an equivalence. This is clear, as the first summand contains all the monomials $x^iy^j$ where $j$ is even, and the second summand those where $j$ is odd. Similarly, we can define maps of $\u$-modules
\[\u\oplus \u[8]\oplus \u[16]\xrightarrow{1\oplus x\oplus x^2} \tmf_7\]
\[\u\oplus \u[8]\oplus \u[12]\oplus \u[16]\oplus\u[20]\oplus\u[28]\xrightarrow{1\oplus x\oplus y \oplus x^2\oplus xy\oplus x^2y}\tmf_{13}\]
at the primes $7$ and $13$, respectively. As in the $p=5$ case, one easily checks these maps are equivalences on homotopy groups. Let us move on to the negative cases now. For $p=11$, we notice that $\pi_\ast \u$ is precisely the summand of $\pi_\ast \tmf_{11}$ supported in nonnegative degrees divisible by $20$. Any potential splitting of $\tmf_{11}$ into sums of $\u$ would have to start by hitting generators in degrees $0$, $8$, $12$, $16$, and $24$. The problem is that we need two summands $\u[24]$ to hit both $y^2$ and $x^3$ in degree $24$ as $\pi_{24}\tmf_{11}=\Z_{11}\{x^3\}\oplus \Z_{11}\{y^2\}$. This means that a potential sum of $u$'s has dimension at least $4$ in degree $64$ as
\[\pi_{64}\u[24]=\pi_{40}\u\simeq \pi_{40}\tmf_{11}\simeq \Z_{11}^2.\] 
This contradicts the fact that the dimension of the $\Z_{11}$-module $\pi_{64}\tmf_{11}$ has dimension $3$. Similar problems happen for primes $p\geq 17$. Indeed, for each of these primes, $\pi_\ast \u$ is the summand of $\pi_\ast \tmf_p$ supported in nonnegative degrees divisible by $2(p-1)$. A potential splitting of $\tmf_p$ into sums of $\u$ would have to hit the two generators in degree $24$, as $2(p-1)\geq 2(16)=32$ is greater than $24$, so $\pi_{24}\u=0$. However, writing $d$ for the dimension
\[d=\dim_{\Z_p} \left(\pi_{2(p-1)}\tmf_p\right) = \dim_{\Z_p}\left(\pi_{2(p-1)}\u \right)\geq 2\]
where the inequality comes from the fact that $2(p-1)\geq 32$, we obtain
\[\dim_{\Z_p}\left(\pi_{2(p-1)+24}\tmf_p\right) = d+1 < 2d = \dim_{\Z_p} \left(\pi_{2(p-1)+24}\left(\u[24]\oplus \u[24]\right)\right).\]
This shows that there can be no splitting of $\tmf_p$ purely in terms of suspensions of $\u$---we will make some suggestions to remedy this in \Cref{previouslydiscussedsplittingconjecturesection}.\\

Onto the periodic case. Consider the basis $\b$ defined as follows: for an even integer $d$, write $\b_d$ for the basis of $\MF_d$ given by
\[\{\Delta^lE_{d'}j^m\}_{m\geq 0}\]
where $d$ is uniquely written as $d=12l+d'$ for $d'$ in the set $\{0,4,6,8,10,14\}$, $j=\frac{x^3}{\Delta}$ is the $j$-invariant, and $E_{d'}$ is the weight $d'$ normalised (meaning with linear term $1$) Eisenstein series which can be summarised by the following formulae:
\[E_0=1\qquad E_4=c_4\qquad E_6=c_6\qquad E_8=c_4^2\qquad E_{10}=c_4c_6\qquad E_{14}=c_4^2c_6\]
Let us write $f_k=\Delta^l E_{k'}$ for the generators of $\MF_k^{\Z_p}$ as a module over $\Z_p[j]\simeq\MF_0^{\Z_p}$. Note these basis elements have some multiplicativity properties which we will implicitly use in what follows:
\[f_{k_1}\cdot f_{12k_2}^r=f_{k_1}\cdot f_{12rk_2}=f_{k_1}\cdot \Delta^{rk_2}=f_{k_1+12rk_2}\]
We now have four cases to consider depending on the remainder of $p$ modulo $12$. Essentially, $f_{p-1}\in \pi_{2(p-1)}\U$ is the first nonzero generator of $\pi_\ast\U$ after $\pi_0\U$. Our splittings of $\TMF_p$ will depend on if $f_{p-1}$ is purely a power of $\Delta$, or a power of $\Delta$ multiplied by $x^2y$, $x$, or $y$. These are precisely the four cases below, respectively.

\paragraph{(The $p\equiv_{12}1$ case)}	Consider the following map of $\U$-modules:
\[\varphi_1\colon \U[2p]\oplus \bigoplus_{\substack{0\leq 2d< p-1 \\ d\neq 1}}\U[4d]\xrightarrow{f_{p}\oplus \bigoplus_d f_{2d}}\TMF_p\]
We claim $\varphi_1$ is an equivalence. First, note the map is injective on homotopy groups, as $\pi_\ast\U$ is concentrated in degrees divisible by $2(p-1)$ and each summand in the domain of the map $\varphi_1$ only hits elements in $\MF_\ast^{\Z_p}$ in degrees which are pairwise distinct modulo $2(p-1)$. In the range $0\leq k\leq p-2$, every $f_k$ is hit by $\varphi_1$ by construction---the only case up for debate is $f_2$, however, $f_{p-1}=\Delta^{\frac{p-1}{12}}$ lies in $\pi_{2(p-1)}\U$ with inverse $\Delta^{\frac{1-p}{12}}$ inside $\pi_{2(1-p)}\U$ and $f_p$ is hit by $\varphi_1$ by construction. We then obtain the following equalities:
\[f_p\cdot f_{p-1}^{-1}=f_p \cdot f_{1-p}=\Delta^{\frac{p-1}{12}-1}x^2y\cdot \Delta^{\frac{1-p}{12}}=\frac{x^2y}{\Delta}=f_2\]
All other $f_k$ are hit for all even $k\in 2\Z$. Indeed, for each such $k$, there is an integer $r$ such that $k+r(p-1)$ lies in the range between $0$ and $p-2$. As $f_{k+r(p-1)}=f_k\cdot f_{p-1}^r$ is hit by $\varphi_1$, and $f_{p-1}^r$ and its inverse lies in $\pi_{2r(p-1)}\U$, we see that the $\pi_\ast\U$-module map induced by $\varphi_1$ hits $f_k$.

\paragraph{(The $p\equiv_{12}11$ case)} Consider the following map of $\U$-modules:
\[\varphi_{11}\colon \bigoplus_{0\leq 2d< p-1}\U[24d]\xrightarrow{\bigoplus f_{12d}=\Delta^d}\TMF_p\]
We claim this map is an equivalence. As in the $p\equiv_{12}1$ case above, we see the induced map on $\pi_\ast$ is injective. To see each $f_k$ in $\MF_\ast^{\Z_p}$ is hit by $\varphi_{11}$, we first note that $f_{6(p-1)}=f_{p-1}^6=\Delta^{\frac{p-1}{2}}$ lies in $\pi_{12(p-1)}\U$, and the above map hits every power of $\Delta$ less than $f_{6(p-1)}$ by construction. In particular, given an even integer $k$, then $f_k$ is hit by $\varphi_{11}$ if $k$ is divisible by $12$. Also, note the following equalities inside $\pi_\ast\U$:
\[f_{p-1}=\Delta^{\frac{p-11}{12}}xy\qquad f_{2(p-1)}=\Delta^{\frac{p-11}{6}+1}x^2\]
\[f_{3(p-1)}=\Delta^{\frac{p-11}{4}+2}y\qquad f_{4(p-1)}=\Delta^{\frac{p-11}{3}+3}x\]
\[f_{5(p-1)}=\Delta^{5\frac{p-11}{12}+3}x^2y\]
If $f_k$ is of the form $\Delta^l E_{k'}$ for $k$ not divisible by $12$, then the equations above show there exists an integer $r$ and an $i$ in the range $1\leq i\leq 5$ such that $f_k=\Delta^r f_{i(p-1)}$ simply because this range of $f_{i(p-1)}$ contain the five remaining possible $E_{k'}$.\\

The following two cases are a mixture of the previous two---let us only detail the first.

\paragraph{(The $p\equiv_{12}5$ case)} Consider the following map of $\U$-modules:
\[\varphi_{5}\colon \bigoplus_{0\leq 2d< p-1}\U[12d]\xrightarrow{\bigoplus f_{6d}}\TMF_p\]
As previously discussed, the induced map on homotopy groups is injective, so it suffices to see $\varphi_5$ hits all the generators of $\MF_\ast^{\Z_p}$. By inspection, we see that $\varphi_5$ hits all $f_k$ of the form $\Delta^i$ and $\Delta^iy$ for all $0\leq i\leq \frac{p-5}{4}$. Moreover, note the following equalities in $\pi_\ast\U$:
\[f_{p-1}=\Delta^{\frac{p-5}{12}}x\qquad f_{2(p-1)}=\Delta^{\frac{p-5}{6}}x^2\qquad f_{3(p-1)}=\Delta^{\frac{p-5}{4}+1}\]
It follows that every $f_k$ of the form $\Delta^i$ and $\Delta^iy$ is hit by $\varphi_5$, for all integers $i$ now. As in the $p\equiv_{12} 11$ case above, the $f_k$'s of the form $\Delta^l x$, $\Delta^lx^2$, $\Delta^l xy$, and $\Delta^lx^2y$, are then hit by $\varphi_5$ as every one of these $E_{k'}$'s is a product of elements in the image-of-$\varphi_5$ by construction or in $\pi_\ast\U$. This shows $\varphi_5$ is an equivalence of $\U$-modules.

\paragraph{(The $p\equiv_{12}7$ case)} The map of $\U$-modules
\[\varphi_{7}\colon \bigoplus_{0\leq 2d< p-1}\U[8d]\xrightarrow{\bigoplus f_{4d}}\TMF_p\]
is an equivalence by an analogous argument to the previous cases.
\end{proof}


\subsection{A conjecture regarding cofibre sequences with $\u$ and $\tmf_p$}\label{previouslydiscussedsplittingconjecturesection}

The fact that $\tmf_p$ is \textbf{not} a quasi-free $\u$-module for primes $p=11$ and $p\geq 17$ seems to be salvageable.

\begin{conjecture}\label{secondpslittingtheorem}
For primes $p\geq 17$ and $p=11$, there exists a cofibre sequence of the following form:
\[\bigoplus_{0\leq 2k<p-1}\u[?]\xrightarrow{\varphi_{\bar{p}}} \tmf_p\to \bigoplus \ell[?]\]
\end{conjecture}

The only real mathematical hurdle left in proving the above conjecture seems to be a combinatorial argument involving the known dimensions of spaces of modular forms of a fixed weight. Let us now see the example for $p=11$ in more detail, and quote the results for $p= 17,19,23$, and $37$.\\

Fix $p=11$ and recall we have the following commutative diagram of $\E_\infty$-rings, a consequence of \Cref{thaversiontwo}:
\[\begin{tikzcd}
{\u}\ar[r]\ar[d]	&	{\tmf_{11}}\ar[d]	\\
{\ell}\ar[r]		&	{\ku_{11}}
\end{tikzcd}\]
Consider the map of $u$-modules
\[y^{10}\colon \u[120]\to \u\]
and its cofibre, which we write as $\u/y^{10}$. By inspecting the homotopy groups of $\u/y^{10}$, one will find they look just like those of the $\u$-module
\[\bigoplus \ell=\ell\oplus \ell[40]\oplus \ell[60]\oplus \ell[80]\oplus \ell[120]\]
viewed with basis $1,x^5,y^5,x^{10}$, and $x^{15}$. To prove that these $u$-modules $\u/y^{10}$ and $\bigoplus \ell$ are equivalent, and more importantly, to later obtain a morphism of $u$-modules from $\ell$ to a quotient of $\tmf_{11}$, consider the cohomological Ext-spectral sequence
\[E_2^{s,t}\simeq \Ext_{\u^\ast}^{s,t}(\pi_{-\ast}M, \pi_{-\ast}N)\Longrightarrow \pi_{-s-t}F_\u(M,N)\]
for any pair of $\u$-modules $M$ and $N$. Setting $M=\u/y^{10}$, the short exact sequence defining $\pi_\ast \u/y^{10}$ shows it has projective dimension 1 as a $\pi_\ast \u$-module, meaning the above spectral sequence is supported in $s=0,1$ and immediately collapses. This degeneration yields a surjection of groups
\[\pi_0 F_\u(\u/y^{10}, N)\to \Ext_{\u^\ast}^{0,0}(\pi_{-\ast}\u/y^{10}, N).\]
Setting $N=\bigoplus \ell$, we lift the desired isomorphism of $\pi_\ast \u$-modules to an equivalence of $\u$-modules
\[\u/y^{10}\simeq \bigoplus \ell.\]
The $\u$-module $\ell$ then naturally maps into $\u/y^{10}$ as the first summand of $\bigoplus \ell$, and with this inclusion we will study a quotient of $\tmf_p$. Consider the following map of $\u$-modules:
\[\varphi_{11}\colon \bigoplus_{d=0}^4 \u[24d]\xrightarrow{1\oplus y^2\oplus y^4\oplus y^6\oplus y^8} \tmf_{11}\]
Write $\tmf_{11}/\varphi$ for the cofibre of this map. Consider the map of $\u$-modules
\[x\colon\u[8]\to \tmf_{11}\]
defined by $x\in \pi_8 \tmf_p$ and the following diagram of $\u$-modules:
\[\begin{tikzcd}
{\u[128]}\ar[r, "y^{10}"]	&	{\u[8]}\ar[d, "{x}"]\ar[r]	&	{\u/y^{10}[8]}\ar[d, dashed]	\\
					&	{\tmf_{11}}\ar[r]			&	{\tmf_{11}/\varphi}
\end{tikzcd}\]
The composite $\u[128]\to \tmf_{11}/\varphi$ vanishes. Indeed, this map of $\u$-modules is represented by the class $xy^{10}$ in $\pi_{128}\tmf_{11}/\varphi$ and $y^{10}=0\in \pi_{120}\tmf_{11}/\varphi$ by the construction of $\varphi_{11}$. Hence, we obtain a map $\u/y^{10}[8]\to \tmf_{11}/\varphi$ which induces multiplication by $x$ on homotopy groups. Precomposing this map with the inclusion $\ell\to \bigoplus \ell$ and the equivalence $\u/y^{10}\simeq \bigoplus \ell$, we obtain the map of $\u$-modules
\[i_x\colon\ell[8]\to \tmf_{11}/\varphi\]
whose effect on homotopy groups is given by multiplication by $x\in \pi_{8}\tmf_{11}$. Replacing $x$ with a class $z\in \pi_{|z|}\tmf_{11}$ in the set
\[Z_{11}=\{y,x^2, x^3, x^4, y^3, x^6, x^7, x^9, x^{12}\}\]
one can repeat the above process, which yields maps of $u$-modules $i_z\colon\ell[|z|]\to \tmf_{11}/\varphi$. These morphisms sum to give the following map of $\u$-modules:
\[i_{11}\colon \bigoplus_{z\in Z_{11}}\ell[|z|]\to \tmf_{11}/\varphi\]
It is now a purely combinatorial exercise to check this is an equivalence. Altogether, this yields the following cofibre sequence of $\u$-modules:
\[\bigoplus_{d=0}^4\u[24d]\xrightarrow{\varphi_{11}} \tmf_{11}\to \bigoplus_{z\in Z_{11}} \ell[|z|]\]
Other examples validating \Cref{secondpslittingtheorem} are the following cofibre sequences:
\[\bigoplus_{d=0}^7 \u[12d]\xrightarrow{\bigoplus y^d}\tmf_{17}\to \bigoplus_{12} \ell[?]\]
\[ \bigoplus_{d=0}^8 \u[8d]\xrightarrow{\bigoplus x^d} \tmf_{19}\to \bigoplus_{9} \ell[?]\]
\[ \bigoplus_{d=0}^{10} \u[24d]\to \tmf_{23}\xrightarrow{\bigoplus \Delta^d} \bigoplus_{55} \ell[?]\]
\[u\oplus \bigoplus_{d=0}^{16} \u[4d+8]\oplus \u[76]\xrightarrow{1\oplus \bigoplus f_{2d+4}\oplus f_{38}} \tmf_{37}\to  \bigoplus_{18} \ell[?]\]
The question marks above signify our lack of understanding of the pattern behind the types of shifts of $\ell$ that occur, although everything above seems to only truly depend upon the residue of the prime modulo $12$.


\subsection{Connective height $2$ image-of-$J$ spectra and \Cref{the}}\label{imageofjapplication}

A classical construction in homotopy theory is that of the \emph{connective image-of-$J$} spectrum $\j$, at the prime $2$ for this exposition, defined by the following cofibre sequence of spectra:
\begin{equation}\label{littlecofibre}\j\to \ko_2\xrightarrow{\psi^3-1}\tau_{\geq 4}\ko_2\end{equation}
This is to be thought of as a connective approximation to the $\Z_2^\times$-fixed points of $\ku_2$, or the $\Z_2^\times/\{\pm 1\}$-fixed points of $\ko$, as $3$ generates $\Z_2^\times$. One constructs the above map by first considering the morphisms of spectra
\begin{equation}\label{compositionofmaps} \ko_2\xrightarrow{\psi^3-1} \ko_2 \to \tau_{\leq 3} \ko_2 \xrightarrow{\simeq} \tau_{\leq 2} \ko_2.\end{equation}
As truncation is a left adjoint, this map is adjoint to $\tau_{\leq 2}\ko_2 \to \tau_{\leq 2}\ko_2$. There is a natural equivalence $ \tau_{\leq 2} \Sph \simeq \tau_{\leq 2}\ko_2$, as $\ko_2$ detects $\eta$ and $\eta^2$, so our desired map is adjoint to $\Sph \to \ko_2 \to \tau_{\leq 2} \ko_2$. This map is zero in degree $0$, hence our original map (\ref{compositionofmaps}) vanishes and hence $\psi^3-1$ factors through $\tau_{\geq 4} \ko_2$ as desired. One can also make an argument using singular cohomology and a Postnikov tower; synthetic versions of this appear in \cite[\textsection4.1]{syntheticj}.\\

To see that $\j$ has a canonical choice of $\E_\infty$-structure, we can also define this spectrum as follows: first, write $\ko_2^{\psi^3}$ for the equaliser of $\psi^3$ and the identity in the category of $\E_\infty$-rings. This is almost $\j$, but there are some stray factors in low degrees, which we deal with by defining $\j$ using the pullback
\begin{equation}\label{heightonepullback}\begin{tikzcd}
{\j}\ar[r]\ar[d]		&	{\ko_2^{\psi^3}}\ar[d]	\\
{\tau_{\leq 2}\Sph}\ar[r]	&	{\tau_{\leq 2}\ko_2^{\psi^3}}
\end{tikzcd}\end{equation}
again in the $\infty$-category of $\E_\infty$-rings---one can check these two definitions match. The practicality of the spectrum $\j$ comes from the fact that the unit map $\Sph_2\to \j$ is split surjective on homotopy groups and detects the $2$-primary \emph{image of the $J$-homomorphism} as well as the Hurewicz image of $\ko$ inside $\pi_\ast\Sph$; this situation is described in \cite{syntheticj} and a simplified proof given too.\\

Here we are interested in defining a height $2$ analogue of the above construction. To this end, we will use the Adams operations of \Cref{thaversiontwo} and adapt (\ref{heightonepullback}) to this $\tmf$-situation. Recall the $\E_\infty$-rings $\u$ of \Cref{heighttwoadamssummands} have Adams operations $\psi^k$ for each $k\in \Z_p^\times$ by \Cref{thaversiontwo}. 

\begin{mydef}\label{higherimageofjstuff}
For any prime $p$, write $g$ for a generator of $\Z_p^\times/F$ where $F$ is the maximal finite subgroup of $\Z_p^\times$. Write $\u^{\psi^g}$ as the equaliser in the $\infty$-category of $\E_\infty$-rings of $\psi^g\colon \u\to \u$ and the identity. At the prime $p=2$, define $\s$ as the $\E_\infty$-ring in the Cartesian square of $\E_\infty$-rings
\[\begin{tikzcd}
{\s}\ar[r]\ar[d]		&	{\u^{\psi^g}=\tmf_2^{\psi^g}}\ar[d]	\\
{\tau_{\leq 6}\Sph_2}\ar[r]&	{\tau_{\leq 6}\u^{\psi^g}.}
\end{tikzcd}\]
For odd primes $p$, let $\s$ be the $\E_\infty$-rings defined by the following Cartesian square:
\[\begin{tikzcd}
{\s}\ar[r]\ar[d]			&	{\u^{\psi^g}}\ar[d]	\\
{\tau_{\leq 2p-3}\Sph_p}\ar[r]	&	{\tau_{\leq 2p-3}\u^{\psi^g}}
\end{tikzcd}\]
\end{mydef}

Notice that $\u^{\psi^g}$ is not connective---the element $1\in \pi_0\u$ contributes to a torsion-free generator $\partial(1)=\zeta\in \pi_{-1}\u^{\psi^g}$. In contrast, the $\E_\infty$-ring $\s$ is connective by construction. The above definition also removes elements of the form $\partial x$ from the homotopy groups of $\u^{\psi^g}$, which leads to the map $\Sph_p\to \s$ being reasonably connective.\\

For example, if $p=2$ and $g=3$, then the map $\Sph_2\to \s$ is an isomorphism on $\pi_d$ for $0\leq d\leq 6$. To see this, let us first calculate some homotopy groups of $\tmf_2^{\psi^3}$. From our knowledge of the Hurewicz image of $\tmf$ (\cite{hurewicztmf} and \cite[\textsection11.11]{brunerrognes}) and the height $1$ image-of-$J$ (\cite[Pr.1.5.22]{greenbook}), we can calculate the following homotopy groups:
\begin{center}
\begin{tabular}{|c|c|c|c|c|c|c|c|c|c|c|c|} 
 \hline
$d$	&	$-1$&	$0$&	$1$&	$2$&	$3$&	$4$&	$5$&	$6$ & $7$ & $8$\\
 \hline
$\text{gen. of }\pi_d \tmf_2^{\psi^3}$	&	$\partial 1$&	$1,\partial\eta$&	$\eta,\partial\eta^2$&	$\eta^2,\partial\nu$&		$\nu$&		$0$&		$\partial\nu^2$&		$\nu^2$	&	$\partial\epsilon,\sigma$ & $\partial c_4\eta,\epsilon$\\
 \hline
\end{tabular}
\end{center}
The elements denoted by $\partial x$ above come from the image of the boundary map $\partial\colon \tmf_2\to \tmf_2^{\psi^3}[1]$ and $\sigma$ is detected by $\partial c_4$, a consequence of the classical height $1$ image-of-$J$ calculation and the fact that the map $\tmf_2\to \ko_2$ commutes with $\psi^3$. From this description, it is clear that $\Sph_2\to \s$ induces an isomorphism on $\pi_d$ for $0\leq d\leq 6$. It is also clear this map is not an isomorphism on $\pi_7$, as $\pi_7\s$ contains $\partial\epsilon$, which does not exist in $\pi_7\Sph_2$. The element $\sigma = \partial c_4$ is the first class detected by $\s$ that is not detected by $\tmf_2$.

\begin{remark}
Some of the utility of $\s$ lie in their relationship to a cofibre sequence akin to (\ref{littlecofibre}). The map $\psi^3-1\colon \tmf_2\to \tmf_2$ factors through $\tau_{\geq 1}\tmf_2$ as $\psi^3$ preserves the unit and $\tau_{\leq 6} \Sph \simeq \tau_{\leq 6} \tmf$, or one can argue with singular cohomology again. The same works at the prime $p=3$ as well---in \cite{heighttwojatthree} we work with a slight variant of $\s$ at the prime $3$.
\end{remark}

We will not pursue similar conclusions at primes $p\geq 5$. Notice that for primes $p\geq 13$ there cannot exist such a simple relationship between cofibres of $\psi^g-1$ and $\s$. Indeed, for such primes, we see that the fibre of a hypothetical map $\psi^g-1\colon \u\to \tau_{\geq 2p-2}$ would have $\pi_{2p-3}$ be a direct sum of $\F_p$'s, one copy for each generator of the ring of modular forms of weight $p-1$. Conversely, by definition, $\pi_{2p-3}\s$ is always $\F_p$, detecting exactly $\al_1$.\\

The following is \Cref{the} and is a complete formality from \Cref{thaversiontwo} and \Cref{higherimageofjstuff}.

\begin{theorem}\label{imageofjdetectionheighttwobody}
Let $p$ be a prime. Then under the unit map $\Sph_p\to \s$, all of the elements in $\pi_\ast \Sph_p$ in the $p$-primary image-of-$J$ and those elements detected by $\Sph\to \tmf_p$ have nontrivial image in $\pi_\ast \s$.
\end{theorem}

In other words, $\pi_\ast\s$ at least detects the Hurewicz image of $\tmf$, which includes the Hurwicz image of $\ko$, and the image of the $J$-homomorphism.

\begin{proof}
Recall that for odd primes, $\j$ is defined either as the fibre of $\psi^g-1\colon \ell\to \ell[2p-2]$ or using the Cartesian diagram of $\E_\infty$-rings
\[\begin{tikzcd}
{\j}\ar[r]\ar[d]			&	{\ell^{\psi^g}}\ar[d]	\\
{\tau_{\leq 2p-3}\Sph}\ar[r]	&	{\tau_{\leq 2p-3}\ell^{\psi^g}}
\end{tikzcd}\]
where $\ell^{\psi^g}$ is the equaliser of $\psi^g$ and the identity. The relationship between $\j$ and the image of the $J$-homomorphism is discussed in \cite[Th.1.5.19 \& Pr.1.5.22]{greenbook} for odd primes and $p=2$, respectively; also see \cite[Th.A]{syntheticj} for a modern discussion and proof. Using this second definition of $\j$, including (\ref{heightonepullback}) at the prime $p=2$, and the fact that $\tmf_p\to \ko_p$ is $\F_p^\times$-equivariant and commutes with $\psi^g$, we naturally obtain a morphism $\s\to \j$ which factors the unit $\Sph_p\to \j$. For our fixed prime $p$, this unit detects the image-of-$J$ in $\pi_\ast\Sph$ and factors through $\s$, so $\s$ also detects the image of the $p$-primary $J$-homomorphism. Similarly, we obtain the detection statement for elements in $\pi_\ast\Sph_p$ detected by $\tmf_p$ as $\Sph_p\to \s$ factors the unit $\Sph_p\to \tmf_p$.
\end{proof}

It would be interesting to know how much closer $\s$ is to $\Sph$ than simply a combination of $\j$ and $\tmf_p$. At primes $p\geq 5$, the spectrum $\s$ does not detect much more than the image-of-$J$. Indeed, a modified Adams--Novikov spectral sequence for $\u^{\psi^g}$ in this case, is concentrated in filtrations $0$ and $1$. However, at the primes $p=2$ and $p=3$ we can still ask the following:

\begin{question}
Are their classes in $\pi_\ast\Sph$ which are detected by $\Sph\to \s$, but which map to zero in $\pi_\ast\j$ and $\pi_\ast\tmf$?
\end{question}

Alternatively, we can ask about the connectivity of $\Sph_p\to \s$. At the prime $p=3$, the author and Christian Carrick show that $\s$ has a rich Hurewicz image; see \cite{heighttwojatthree}. In current work-in-progress, we also explore the connection between the Hurewicz image of $\s$ as the elements constructed in \cite{afewinfinitefamilies}. In particular, the answer to the above question in this case is an emphatic ``yes!''.

\addcontentsline{toc}{section}{References}


\bibliography{/Users/jackdavies/Dropbox/Work/references} 

\bibliographystyle{alpha}


\end{document}